\documentclass[11pt]{article}
\usepackage{amsthm,amssymb,amsmath,latexsym}
\usepackage{graphs}
\theoremstyle{definition}
\newtheorem{defin}{Definition}[section]
\theoremstyle{remark}

\newtheorem*{remar}{Remark}
\theoremstyle{plain}
\newtheorem{thm}[defin]{Theorem}
\newtheorem{prop}[defin]{Proposition}
\newtheorem{lemm}[defin]{Lemma}
\newtheorem{corol}[defin]{Corollary}

\numberwithin{equation}{section}

\begin{document}

\title{Homeomorphic Measures on Stationary Bratteli Diagrams }

\author{ S.~Bezuglyi and O.~Karpel\footnote{O. Karpel was supported in part by the Akhiezer fund.}\\
Institute for Low Temperature Physics, \\
47 Lenin Avenue, 61103 Kharkov, Ukraine \\
(e-mail: bezuglyi@ilt.kharkov.ua, helen.karpel@gmail.com)}

\date{}
\maketitle

\begin{abstract} We study the set $\mathcal S$ of ergodic probability  Borel measures on stationary non-simple Bratteli diagrams which are invariant with respect to the tail equivalence relation. Equivalently, the set $\mathcal S$ is formed by ergodic probability measures invariant with respect to aperiodic substitution dynamical systems. The paper is devoted to the classification of  measures $\mu$ from $\mathcal S$ with respect to a homeomorphism. The properties of these measures related to the clopen values set $S(\mu)$ are studied. It is shown that for every measure $\mu \in \mathcal S$ there exists a subgroup $G \subset \mathbb R$ such that $S(\mu) = G\cap [0,1]$, i.e. $S(\mu)$ is group-like. A criterion of goodness is proved for such measures. Based on this result, we classify the measures  from $\mathcal S$ up to a homeomorphism. It is proved that for every good measure $\mu\in \mathcal S$ there exist countably many measures $\{\mu_i\}_{i\in \mathbb N}\subset \mathcal S$ such that the measures $\mu$ and $\mu_i$ are homeomorphic but the tail equivalence relations on the corresponding Bratteli diagrams are not orbit equivalent.

\end{abstract}

\section{Introduction}

In this paper, we are interested in the problem of classification of Borel probability measures  on a Cantor set with respect to  a homeomorphism. Two probability measures $\mu$ and $\nu$ defined on Borel subsets of a topological space $X$ are called \textit{homeomorphic} or \textit{topologically equivalent} if there exists a self-homeomorphism $h$ of  $X$  such that $\mu = \nu\circ h$, i.e. $\mu(E) = \nu(h(E))$ for every Borel subset $E$ of $X$. In such a way, the set of all Borel probability measures on $X$ is partitioned into equivalence classes. One may be interested in the structure of the equivalence relation defined by the classes of homeomorphic measures or in the study of a certain equivalence class.

The topological properties of the space $X$ are important for the classification of measures up to a homeomorphism. For instance, the following theorem proved by Oxtoby and Ulam \cite{Oxt-Ul} holds: a non-atomic Borel probability  measure $\mu$ on the finite-dimensional cube $[0, 1]^{n}$ is homeomorphic to the Lebesgue measure if and only if every nonempty open set has a positive measure (in other words, $\mu$ is full) and the boundary of the cube has measure 0. Later, Oxtoby and Prasad extended this result to the Hilbert cube $[0, 1]^{\mathbb N}$. Similar results were also obtained for various manifolds (see the book by Alpern and Prasad \cite{Alp-Pr} for the details).

The current work has two sources. The first one is the article \cite{S.B.} where an explicit description of all ergodic (finite and infinite)  measures on  stationary Bratteli diagrams was found. The second one is a series of papers by Akin, Austin, Dougherty, Mauldin, Yingst (\cite{Akin3, Austin, D-M-Y, Yingst}) where Borel probability measures on zero-dimensional compact perfect metric spaces (Cantor sets) were extensively studied. In those papers, the major results were focused on the classification of Bernoulli measures up to a homeomorphism continuing the preceding investigations (see \cite{Navarro, Navarro-Oxtoby}).  We should mention that it was Akin who initiated a systematic study  of  homeomorphic measures on a Cantor space \cite{Akin1, Akin2}.   It turns out that the situation in this case is much more difficult than for connected spaces. Though there is, up to a homeomorphism, only one Cantor set it is not hard to construct full non-atomic  measures on the Cantor set which are not topologically equivalent (see \cite{Akin1} where more impressive results were proved). The following question naturally  arises: find necessary and sufficient  conditions under which measures on a Cantor space $X$ are homeomorphic. It was  noted in  \cite{Akin1}  that  there exist continuum  classes of equivalent full non-atomic probability measures on a Cantor set. This fact is based on the existence of a countable base of clopen subsets of a  Cantor set. Akin \cite{Akin1} defined the \textit{clopen values set} $S(\mu)$ as the set of values of measure $\mu$ on all clopen subsets of $X$. The set $S(\mu)$ is a countable dense subset of the unit interval, and this set provides an invariant for topologically equivalent measures, although it is not a complete invariant, in general. But for the class of the so called \textit{good} measures, $S(\mu)$ \textit{is} a complete invariant. By definition, a  full non-atomic probability measure $\mu$ is good if whenever $U$, $V$ are clopen sets with $\mu(U) < \mu(V)$, there exists a clopen subset $W$ of $V$ such that $\mu(W) = \mu(U)$. It turns out that such measures are exactly invariant measures of uniquely ergodic minimal homeomorphisms of Cantor sets (see \cite{Akin2}, \cite{GW}).

It this paper, we consider stationary (non-simple) Bratteli diagrams and ergodic probability measures on their path spaces invariant with respect to the cofinal (tail) equivalence relation. It follows  from \cite{S.B.} that, for every such measure $\mu$, the key invariant $S(\mu)$ can be easily computed in terms of eigenvector entries and eigenvalues of the corresponding incidence matrices. This allows us to answer the questions about properties of $\mu$ and $S(\mu)$ and construct homeomorphic measures.

The proved results and organization of the paper are the following. In Section \ref{AnSurv} we collect the definitions and statements about measures on a Cantor set that are used in the paper. Since we would like to make the paper self-contained, we include  the main results from \cite{Akin2}, \cite{Akin3}, \cite{S.B.} in this section. We also discuss the main notions related to Bratteli diagrams. In contrast to the case of full measures on a Cantor set, we have to deal with several singular ergodic invariant measures on a path space of a stationary Bratteli diagram. As shown in \cite{S.B.}, the support of every such measure is a closed (Cantor) subset and any open subset of the support has positive measure. This allows us to use the machinery  developed for full measures.
In Section \ref{section3}, we first study the structure of the clopen values set $S(\mu)$ for any ergodic probability $\mathcal R$-invariant measure $\mu$ on a stationary Bratteli diagram $B$ and  prove that this set is group-like, i.e., $S(\mu) = G\cap [0,1]$ for some additive subgroup  $G$ of $\mathbb R$.
It is worth to mention that in the main result of Section \ref{section3} we  consider two cases:  (i) $S(\mu)$ is a subset of $\mathbb Q$, and (ii) $S(\mu)\cap (\mathbb R \setminus \mathbb Q) \neq \emptyset$. The first case is relatively simple. But in the second case we have to use some methods of linear algebra and matrix theory. As proved in \cite{S.B.}, every ergodic finite invariant measure is completely determined in terms of eigenvalues and eigenvectors of the matrix that defines the Bratteli diagram $B$. In this case, the eigenvector entries and eigenvalues admit their representations as vectors with rational entries.  In the same section, we prove an easy checkable criterion for a measure from $\mathcal S$ to be good. In Section \ref{section4},  we apply the found criterion answering the following question. Given a good measure $\mu$ on a stationary Bratteli diagram, how many measures from $\mathcal S$ are homeomorphic to $\mu$? Is this class infinite? It is proved that
there exist stationary Bratteli diagrams $\{B_i\}_{i=0}^\infty$ and  good ergodic $\mathcal{R}_i$-invariant probability measures $\mu_i$ on $B_i$ such that each measure $\mu_i$ is homeomorphic to the given measure $\mu$, but the dynamical systems $(B_i, \mathcal R_i)$, $(B_j, \mathcal R_j)$ are topologically orbit equivalent if and only if $i = j$. The last section contains several examples that illustrate the results proved in the preceding sections. Namely, we give a class of stationary non-simple Bratteli diagrams such that in the simplex of probability $\mathcal R$-invariant measures only ergodic measures are good. Another example of a set Bratteli diagrams contains an explicit description of all good measures. It is also shown that given a measure $\nu$ on a stationary Bratteli diagram there exists a good measure $\mu$ such that $S(\mu) = S(\nu)$.


\section{Preliminaries: Good Measures and Stationary Bratteli Diagrams}\label{AnSurv}

In this section, we collect some necessary definitions and results that are used throughout the paper. We do this for the reader's convenience. Unless stated otherwise, all measures considered in the paper are \textit{Borel probability non-atomic} measures and all Bratteli diagrams are assumed to be \textit{stationary}.
\medskip

{\bf 2.1. Good measures}. For a measure $\mu$ on a Cantor  space $X$, define the
\textit{clopen values set:}
$$
S(\mu) = \{\mu(U):\,U\mbox{ is clopen in } X\}.
$$
For each probability measure $\mu$ on $X$, the set $S(\mu)$ is a dense subset of the unit interval containing  $0$ and $1$ \cite{Akin1}.

Let $X_{1}$, $X_{2}$ be two Cantor sets, $h \colon X_{1} \rightarrow X_{2}$ a continuous map, and $\mu_{1}$ a  measure on $X_{1}$. Then the \textit{image measure} $h_{*}\mu_{1}$ on $X_{2}$ is defined by
$$
h_{*}\mu_{1}(B) = \mu_{1}(h^{-1}(B))
$$
for all Borel subsets $B$ of $X_{2}$. It is said that the measures $\mu_{1}$ on $X_{1}$ and $\mu_{2}$ on $X_{2}$ are \textit{homeomorphic} if there exists a homeomorphism $h \colon X_{1} \rightarrow X_{2}$ such that $h_{*}\mu_{1} = \mu_{2}$. Clearly, $S(\mu_{1}) = S(\mu_{2})$ for any homeomorphic measures $\mu_1$ and $\mu_2$.

A measure  $\mu$ on a Cantor set $X$ is called \textit{full} if $\mu(V) > 0$ for any non-empty clopen subset $V$ of $X$. If $\mu(V) > 0$, then one can  define the relative measure $\mu_{V}$ on $V$ setting
$$
\mu_{V}(A) =\frac{\mu(A\cap V)}{\mu(V)}
$$
where $A$ is a Borel subset of $X$.

We recall below the  definitions of \textit{good, refinable, and weakly refinable measures} which are based on some natural properties of measures on a Cantor set. We follow here the papers  \cite{Akin2}, \cite{Akin3}, and  \cite{D-M-Y}.

A \textit{partition basis} $\mathcal{B}$ for a Cantor set $X$ is a collection of clopen subsets of $X$ such that every non-empty clopen subset of $X$ can be partitioned by elements of $\mathcal{B}$. A partition basis is a basis for the topology but not every basis is a partition basis.

\begin{defin}
Let $\mu$ be a full measure on a Cantor set $X$.

  (1) A clopen subset $V$ of $X$ is called \textit{good} for $\mu$ (or just good when the measure is understood) if for every clopen subset $U$ of $X$ with $\mu(U) < \mu(V)$, there exists a clopen set $W$ such that $W \subset V$ and $\mu(W) = \mu(U)$.
A measure $\mu$ is called \textit{good} if every clopen subset of $X$ is good for $\mu$.

 (2) A clopen subset $U$ of $X$ is called \textit{refinable} for $\mu$  if $\alpha_{1}, \ldots ,\alpha_{k} \in S(\mu)$ with $\alpha_{1}+ \ldots + \alpha_{k} = \mu(U)$ implies that there exists a clopen partition $\{U_{1}, \ldots , U_{k}\}$ of $U$ with $\mu(U_{i}) = \alpha_{i}$ for $i = 1,\ldots ,k$.
A measure $\mu$ is called \textit{refinable} if every clopen subset is refinable.

 (3) A measure $\mu$ is called \textit{weakly refinable} if there exists a partition basis $\mathcal{B}$ for $X$ with $X \in \mathcal{B}$ consisting of refinable clopen subsets.

 (4) A non-empty clopen subset $U$ of $X$ is called a \textit{clopen set of $\mu$ type} when $S(\mu_{U}) = S(\mu)$.

 (5) A measure $\mu$ is called a measure of \textit{Bernoulli type} when there is a partition basis $\mathcal{B}$ for $X$ consisting of clopen sets of $\mu$ type.

 (6) It is said that a measure $\mu$ on a Cantor set $X$ satisfies the \textit{Quotient Condition} when every non-empty clopen subset $U$ of $X$ is of $\mu$ type.

 (7) A subset $S$ of the unit interval $I = [0, 1]$ is called \textit{group-like (ring-like, field-like)} if $S = G \bigcap I$ where  $G$ is an additive subgroup (subring, subfield) of $\mathbb{R}$.
\end{defin}

We note that ``goodness'' $\Longrightarrow$ ``refinability''  $\Longrightarrow$ ``weak refinability''. One can  find refinable but not good measures. To the best of our knowledge, it is an open question whether the notions of refinability and weak refinability are equivalent.

It can be easily verified that for a countable subset $S$ of the unit interval with $0, 1 \in S$ the set $S$ is group-like (ring-like, field-like) if and only if $S + \mathbb{Z}$ is a subgroup (subring, subfield) of $\mathbb{R}$. Clearly, $S$ is ring-like if and only if $S$ is group-like and multiplicative. In fact, the following lemma proved by Akin \cite{Akin2} holds.

\begin{lemm}\label{AkinKritGrouplike}
Let $S$ be a subset of $[0,1]$ with $0,1 \in S$. Let $G(S)$ be the additive group of $\mathbb{R}$ generated by $S$. The following conditions on $S$ are equivalent:

(1) $S$ is group-like;

(2) $S + \mathbb{Z} = G(S)$;

(3) $S + \mathbb{Z}$ is an additive subgroup of $\mathbb{R}$;

(4) $\alpha, \beta \in S$ and $\alpha \leq \beta$ imply that $\beta - \alpha \in S$.

\noindent
If $S$ is group-like and $G$ is an additive subgroup of $\mathbb{R}$, then
$$
S = G \bigcap [0,1] \Leftrightarrow S + \mathbb{Z} = G.
$$

\end{lemm}

In the next theorem we collect the following results on good measures.

\begin{thm} \label{goodmeasure}
Let $\mu$ be a  measure on a Cantor set $X$.

(1) If $\mu$ is good and a measure $\nu$ is homeomorphic to $\mu$, then $\nu$ is good.

(2) In order that the measure $\mu$ on $X$ be good, it suffices that there exists a partition basis $\mathcal B$ consisting of clopen sets which are good for $\mu$. In particular, if a clopen set can be partitioned by good clopen sets, then it is itself good \cite{Akin3}.

(3) The measure $\mu$ is good if and only if there is a uniquely ergodic, minimal homeomorphism of the Cantor set for which $\mu$ is the unique invariant measure \cite{Akin2}, \cite{GW}.

(4) The direct product of a finite or infinite sequence of good measures is a good measure \cite{Akin3}.

(5)  If $\mu$ is good, then $S(\mu)$ is group-like. Conversely, if $S$ is a group-like countable dense subset of $[0,1]$, then there is a good measure $\mu$ on $X$ such that $S = S(\mu)$ ($\mu$ is unique up to a homeomorphism).

(6) If $\mu$ is good and $V$ is a non-empty clopen subset of $X$, then $\mu_{V}$ is a good measure on the Cantor set $V$ and therefore $S(\mu_{V})$ is group-like \cite{Akin2}.

(7) \label{gd-ref} The following statements are equivalent \cite{Akin3}:\\
 (i) $\mu$ is a good measure;\\
 (ii) $\mu$ is refinable and $S(\mu)$ is group-like;\\
 (iii) $\mu$ is weakly refinable and $S(\mu)$ is group-like.
\end{thm}

\begin{corol}\label{homeomorphic_measures} If $\mu$ and $\nu$ are good measures
on Cantor sets $X$ and $Y$ and if $S(\mu) = S(\nu)$, then $\mu$ and $\nu$ are  homeomorphic.
\end{corol}

Let $D$ be a countable subset of the unit interval which contains 1. A number $\delta \in [0, 1]$ is called a \textit{divisor} of $D$ if for all $\alpha \in [0, 1]$
$$
\alpha \in D \Longleftrightarrow \alpha\cdot\delta \in D.
$$
The set of all divisors of $D$ is denoted by $Div(D)$. The set $Div(D)$ is multiplicative and $1 \in Div(D)$.

This following theorem focuses on the properties of measures of Bernoulli type:

\begin{thm}\label{BernQuot}
Let $\mu$ be a full  measure on a Cantor set $X$ and let $G$ be the group generated by the clopen values set $S(\mu)$.

(1) If $\mu$ is good then $\mu$ is of Bernoulli type if and only if $G$ is a subring of $\mathbb{R}$ such that every positive element of $G$ is a sum of positive units of $G$.

(2) If $\mu$ is good then every non-empty clopen subset of $X$ is a set of $\mu$ type if and only if $G$ is a subfield of $\mathbb{R}$.

(3) If $\mu$ is of Bernoulli type then the clopen values set $S(\mu)$ is multiplicative and for every non-empty clopen subset $U$ of $X$ we have $S(\mu) \subset S(\mu_{U})$.

(4) If $\mu$ satisfies the Quotient Condition then $\mu$ is a refinable measure of Bernoulli type, and, for every non-empty clopen $U \subset X$, the relative measure $\mu_{U}$ is homeomorphic to $\mu$. The rationals $\mathbb{Q} \bigcap [0,1]$ are contained in $S(\mu)$.

(5) Any two of the following conditions imply the third:\\
(i) $U$ is a clopen subset good for $\mu$;\\
(ii) $U$ is a clopen subset of $\mu$ type;\\
(iii)  $\mu(U) \in Div(S(\mu))$.
\end{thm}

\textbf{2.2. Bratteli diagrams}. We recall here some basic definitions and facts about Bratteli diagrams.  We mainly use the notation and results from \cite{S.B.}.

\begin{defin}
A \textit{Bratteli diagram} is an infinite graph $B = (V,E)$ such that the vertex set $V = \bigcup_{i\geqslant 0}V_{i}$ and the edge set $E = \bigcup_{i\geqslant 1}E_{i}$ are partitioned into disjoint subsets $V_{i}$ and $E_{i}$ such that

(i) $V_{0} = \{v_{0}\}$ is a single point;

(ii) $V_{i}$ and $E_{i}$ are finite sets;

(iii) there exist a range map $r$ and a source map $s$ from $E$ to $V$ such that $r(E_{i}) = V_{i}$, $s(E_{i}) = V_{i-1}$, and $s^{-1}(v)\neq 0$, $r^{-1}(v')\neq 0$ for all $v \in V$ and $v' \in V \setminus V_{0}$.
\end{defin}

The pair $(V_{i}, E_{i})$ or just $V_i$ is called the $i$-th level of the diagram $B$.
A finite or infinite sequence of edges $(e_{i} : e_{i} \in E_{i})$ such that $r(e_{i}) = s(e_{i + 1})$ is called a \textit{finite} or \textit{infinite path}, respectively. For a Bratteli diagram $B$, we denote by $X_{B}$ the set of all infinite paths starting at the vertex $v_{0}$. We endow $X_{B}$ with the topology generated by cylinder sets $U(e_{1}, \ldots ,e_{n}) = \{x \in X_{B} : x_{i} = e_{i}, i = 1, \ldots , n \}$, where $(e_{1}, \ldots ,e_{n})$ is a finite path from $B$.  We consider here such Bratteli diagrams $B$ for which the path space $X_{B}$ is a Cantor set.

Given a Bratteli diagram $B = (V,E)$, define a sequence of incidence matrices  $F_{n} = (f_{vw}^{(n)})$ of $B$:
$$
f_{vw}^{(n)} = |\{e \in E_{n+1} : r(e) = v, s(e) = w\}|
$$
where $v \in V_{n+1}$ and $w \in V_{n}$ and the size of $F_n$ is  $|V_{n+1}|\times |V_{n}|$. Here and thereafter $|\Lambda|$ denotes the cardinality of the set $\Lambda$.

A Bratteli diagram is called \textit{stationary} if $F_{n} = F_{1}$ for every $n \geq 2$.

Observe that every vertex $v \in V$ is connected to $v_{0}$ by a finite path, and the set $E(v_0,v)$ of all such paths is finite. Set $h_{v}^{(n)} = |E(v_{0}, v)|$, $v \in V_{n}$. Then
$$
h^{(n + 1)} = F_{n}h^{(n)}.
$$
where $h^{(n)} = (h_{w}^{(n)})_{w \in V_{n}}$.

For $w \in V_{n}$, the set $E(v_{0}, w)$ defines the clopen subset $X_{w}^{(n)} := \{x = (x_{i}) \in X_{B} : r(x_{n}) = w \}$ of $X_B$. Then  $\{X_{w}^{(n)} : w \in V_{n}\}$ is a clopen partition of $X_{B}$. Analogously, the sets $X_{w}^{(n)}(\overline{e}) := \{x = (x_{i}) \in X_{B} : x_{i} = e_{i}, i = 1,...,n\}$ determine a clopen partition of $X_{w}^{(n)}$ where  $\overline{e} = (e_{1}, \ldots ,e_{n}) \in E(v_{0}, w)$, $n \geq 1$.

\begin{defin}
Let $B = (V,E)$ be a Bratteli diagram. Two infinite paths $x = (x_{i})$ and $y = (y_{i})$ from $X_{B}$ are called \textit{tail equivalent} if there exists $i_{0}$ such that $x_{i} = y_{i}$ for all $i \geq i_{0}$. Denote by $\mathcal{R}$ the tail equivalence relation on $X_{B}$.
\end{defin}

Recall that a Bratteli diagram is called \textit{simple} if the tail equivalence relation $\mathcal R$ is minimal.

We will consider Bratteli diagrams $B$ for which $\mathcal{R}$ is a countable Borel equivalence relation on $X_{B}$. Any two paths $x$, $y$ from $X_{B}$ are $\mathcal{R}$-equivalent if and only if there exists $w \in V$ such that $x \in X_{w}^{(n)}(\overline{e})$ and $y \in X_{w}^{(n)}(\overline{e'})$ for some $\overline{e}, \overline{e'} \in E(v_{0}, w)$.

Recall that a measure $\mu$ on $X_B$ is called $\mathcal{R}$-\textit{invariant} if for any two paths $\overline{e}$ and $\overline{e'}$ from $E(v_{0}, w)$ and any vertex $w$, one has $\mu(X_{w}^{(n)}(\overline{e})) = \mu(X_{w}^{(n)}(\overline{e'}))$.  Then
$$
\mu(X_{w}^{(n)}(\overline{e})) =  \frac{1}{h_{w}^{(n)}}\mu(X_{w}^{(n)}), \ \  \overline{e} \in E(v_{0}, w).
$$

In \cite{S.B.}, all invariant ergodic measures on a stationary Bratteli diagram were  described as follows. It was first shown  that the study of any stationary Bratteli diagram $B = (V,E)$ (with $|V| = K$) can be reduced to the case when the incidence matrix $F$ of size $K\times K$ has the form:

\begin{equation}\label{Frobenius Form}
F =\left(
  \begin{array}{ccccccc}
    F_1 & 0 & \cdots & 0 & 0 & \cdots & 0 \\
    0 & F_2 & \cdots & 0 & 0 & \cdots & 0 \\
    \vdots & \vdots & \ddots & \vdots & \vdots & \cdots& \vdots \\
    0 & 0 & \cdots & F_s & 0 & \cdots & 0 \\
    X_{s+1,1} & X_{s+1,2} & \cdots & X_{s+1,s} & F_{s+1} & \cdots & 0 \\
    \vdots & \vdots & \cdots & \vdots & \vdots & \ddots & \vdots \\
    X_{m,1} & X_{m,2} & \cdots & X_{m,s} & X_{m,s+1} & \cdots & F_m \\
  \end{array}
\right)
\end{equation}

The square non-zero matrices $F_i$, $i=1,...,m$, are irreducible (without loss of generality, one can assume that these matrices are strictly positive). For any $j = s+1,...,m$, at least one of the matrices $X_{j,k}$ is non-zero. The matrices $F_i$ determine the partition of the vertex set $V$ into subsets (classes) $V_i$ of vertices. In their turn, these subsets generate  subdiagrams $B_i$. The non-zero matrices $X_{j,k}$ indicate which subdiagrams  are linked by some edges (or finite paths). Notice that each subdiagram $B_i$, $i=1,...,s$, corresponds to a minimal component of the cofinal equivalence relation $\mathcal R$.

We denote by $F_\alpha$, $\alpha \in \Lambda$, the \textit{non-zero} matrices on the main diagonal in (\ref{Frobenius Form}).
Let $\alpha \geq \beta$. It is said that the class of vertices $\alpha$ {\it has access} to a class $\beta$, in symbols $\alpha \succeq \beta$, if and only if either $\alpha = \beta$ or there is a finite path in the diagram from a vertex which belongs to $\beta$ to a vertex from  $\alpha$. In other words, the matrix $X_{\alpha, \beta}$ is non-zero. A class $\alpha$ is called {\it final (initial)} if there is no class $\beta$ such that $ \alpha \succ   \beta$ ($\beta \succ \alpha$).

Let  $\rho_\alpha$ be the spectral radius of $F_\alpha$. A class $\alpha\in \{1,...,m\}$ is called  {\it distinguished} if $\rho_\alpha > \rho_\beta$ whenever $ \alpha \succ   \beta$. Notice that all classes $\alpha = 1,\ldots,s$ are necessarily distinguished. A real number $\lambda$ is called a {\it
distinguished eigenvalue} if there exists a non-negative eigenvector $x$ with $Fx=\lambda x$.
If $x = (x_1,...,x_K)^T$ is an eigenvector corresponding to a distinguished eigenvalue $\lambda_\alpha$, then $x_i >0$ if and only if $i \in \beta$ and $\alpha \succeq \beta$.

Let  $\lambda_1,...,\lambda_k$ be the distinguished eigenvalues of the matrix $A = F^T$ (we will keep this notation below).
The main result of \cite{S.B.} asserts that there exist exactly $k$ ergodic probability invariant measures defined by  $\lambda_1,...,\lambda_k$. More precisely, fix a distinguished eigenvalue $\lambda$ and let  $x =  (x_1,...,x_K)^T$ be the probability non-negative eigenvector  corresponding to  $\lambda$. Then the ergodic probability  measure $\mu$ defined by $\lambda$ and $x$ satisfies the relation:
\begin{equation}\label{mucount}
\mu(X_{i}^{(n)}(\overline{e})) = \frac{x_i}{\lambda^{n - 1}}
\end{equation}
where $i\in V_n$ and $\overline e$ is a finite path with $s(\overline e) =i$. Therefore, the clopen values set for $\mu$ has has the form:
\begin{equation}\label{formulaSmu}
S(\mu) = \left\{\sum_{i = 1}^{K} k^{(n)}_{i}\frac{x_i}{\lambda^{n - 1}} : 0  \leq k^{(n)}_{i} \leq h^{(n)}_{i}; \; n = 1, 2, \ldots \right\}.
\end{equation}
This relation is of extreme importance for us and will be used throughout the paper.

Let $\lambda_\alpha$ be a distinguished eigenvalue corresponding to the distinguished class $\alpha$.
In the next section, we will use the following asymptotics mentioned in   \cite{S.B.}
\begin{equation} \label{eq-asymp}
(A^n)_{ij} \sim \lambda_\alpha^n,\ \ n\to \infty,\ \ \mbox{for}\ i\in \beta,\ j\in \alpha,\ \mbox{with}\
\alpha \succeq \beta.
\end{equation}
Here $\sim$ means that the ratio tends to a positive constant. On the other hand,
\begin{equation} \label{eq-asymp2}
(A^n)_{ij} = o(\lambda_\alpha^n),\ \ n\to \infty,\ \ \mbox{for}\ j\in
\beta \prec \alpha.
\end{equation}

If $\lambda$ is a non-distinguished Perron-Frobenius eigenvalue for $A$, then the corresponding $\mathcal R$-invariant measure on $X_B$ is infinite \cite{S.B.}. We do not study infinite measures in this paper.
\medskip

\textbf{ 2.3. Measure supports}. Given the diagram $B$ as above, let $Y_\alpha$ be the path space of the Bratteli subdiagram $B_\alpha,\
\alpha \in \Lambda$.  Define $X_\alpha = \mathcal
R(Y_\alpha)$, that is, a path $x\in X_B$ belongs to $X_\alpha$ if it
is $\mathcal R$-equivalent to a path $y\in Y_\alpha$.  We see that
$X_\alpha = Y_\alpha$ if and only if  $\alpha$ is a distinguished class corresponding to a minimal component of $\mathcal
R$. It follows from the structure of the diagram $B$, see (\ref{Frobenius Form}), that $X_\alpha \cap X_\beta =
\emptyset$ for $\alpha \neq \beta$, and  $\{X_\alpha : \alpha \in
\Lambda\}$ is a partition of $X_B$.
It is also easy to see that for any $x\in X_\alpha$ the orbit
$\mathcal R(x)$ is dense in $\bigcup_{\beta\prec\alpha} X_\beta$.

We describe here the support of measure $\mu_\alpha$ defined by
a distinguished eigenvalue $\lambda_\alpha$ and the corresponding eigenvector $x_\alpha$ of $A$. Given such a measure $\mu_\alpha$, we call the measure support, $supp(\mu_\alpha)$, the largest closed (Cantor) subset of $X_B$ such that every open subset of $supp(\mu_\alpha)$ has positive measure. In other words, $\mu_\alpha$ is full on $supp(\mu_\alpha)$.

It is clear that when $\alpha$ is a final class, then the support of $\mu_\alpha$ is a Cantor set $Y_\alpha$.

If $\alpha$ is not a final class, then the measure $\mu_\alpha$ is sitting on $X_\alpha = \mathcal R(Y_\alpha)$. Then $supp(\mu_\alpha)$ is the closure $\overline{X_\alpha}$  of $X_\alpha$. That is, to obtain $\overline{X_\alpha}$,  we need to add to $X_\alpha$ those minimal components of the diagram which are accessible from $\alpha$.

It is obvious that all definitions (clopen values set, good measures,  measures of Bernoulli type,  etc.)  given in subsection \textbf{2.1} are applicable to the measures $\mu_\alpha$ because $\mu_\alpha$ is full on the corresponding Cantor sets $\overline X_\alpha$. In particular, we note that the ergodic measures corresponding to minimal components are automatically good: on a simple stationary Bratteli diagram any Vershik map is minimal and uniquely ergodic.

Finally, we remark that one can extend the mentioned above definitions to non-ergodic $\mathcal R$-invariant measures. Such measures form a simplex whose extreme points are ergodic $\mathcal R$-invariant measures. Therefore, any non-ergodic finite $\mathcal R$-invariant measure is supported on the closure of a finite disjoint union of some  sets $X_\alpha$.


\section{Good Measures on Stationary Bratteli diagrams}\label{section3}

In this section, we study finite ergodic $\mathcal{R}$-invariant measures on stationary Bratteli diagrams. We show that, for any such a measure, the clopen values set $S(\mu)$ is group-like. We also give the necessary and sufficient conditions under which a measure on a stationary Bratteli diagram is good.
\medskip

\textbf{3.1. Group-like clopen values set}. Consider a stationary non-simple Bratteli diagram $B = (V, E)$. Let $F$ be its incidence $K \times K$ matrix and $A = F^T$. Let $\mu$ be the measure defined by a distinguished class of vertices $\alpha$ and $\lambda$  the corresponding distinguished eigenvalue of $A$.
Denote by  $(y_1,...,y_K)^T$  the probability  eigenvector of the matrix $A$ corresponding to $\lambda$. Notice that the vector $(y_1,...,y_K)^T$ may have zero entries. These zero entries are assigned  to the vertices from $B$ that are not accessible from the class $\alpha$. Denote by $(x_1, \ldots, x_n)^T$ the positive vector obtained from $(y_1,...,y_K)^T$ by crossing out zero entries.  We call  $(x_1,...,x_n)^T$ the \textit{reduced vector} corresponding to the measure $\mu$.  Recall that we consider the measure $\mu$ only on its support. This means that we can ignore the part of  $B$ formed by subdiagrams which are not accessible from the  class $\alpha$. Without loss of generality, we can think that the  matrix $A = F^T$ satisfies the condition $Ax = \lambda x$.

Let $H$ be the additive subgroup of $\mathbb{R}$ generated by $\{x_1, \ldots , x_n\}$.

\begin{lemm} \label{lambdaH}
Let $B$, $\mu$, $A$,  $\lambda$, $H$, and $(x_1,...,x_n)^T$ be as above.  Let $G$ be the additive subgroup of $\mathbb{R}$ generated by the clopen values set $S(\mu)$. Then:

1) $\lambda H \subset H$;

2) $G = \bigcup\limits_{N \in \; \mathbb{N}} \frac{1}{\lambda^N}H$;

3) $S(\mu)$ is group-like if and only if $S(\mu) + \mathbb{Z} = \bigcup\limits_{N \in \; \mathbb{N}} \frac{1}{\lambda^N}H$;

4) $\lambda^{-1} \in Div(G)$.

\end{lemm}

\noindent
\textbf{Proof.}
1) Let $A = (a_{ij})_{i,j = 1}^{n}$. Then $\sum_{j=1}^{n} a_{ij}x_j = \lambda x_{i}$, hence  $\lambda x_{i} \in H$ for $i = 1, \ldots, n$. Since $\lambda H$ is generated by
$\lambda x_1,...,\lambda x_n$, we see that $\lambda H \subset H$.

 2) It follows from (\ref{formulaSmu}) that $\frac{1}{\lambda^N}H \subset G$ because $\frac{x_1}{\lambda^{N}}, \ldots, \frac{x_n}{\lambda^{N}} \in G$ for any $N \in \mathbb{N}$. On the other hand, we see that $\bigcup\limits_{N \in \; \mathbb{N}} \frac{1}{\lambda^N}H$ is a group and $S(\mu) \subset \bigcup\limits_{N \in \; \mathbb{N}} \frac{1}{\lambda^N}H$. Hence, $G \subset \bigcup\limits_{N \in \; \mathbb{N}} \frac{1}{\lambda^N}H$.

3) By Lemma \ref{AkinKritGrouplike}, $S(\mu)$ is group-like if and only if $G = S(\mu) + \mathbb{Z}$.

4) It follows from the above results that $\lambda G = G$. \hfill$\blacksquare$

\begin{remar}
Since $\sum_{k=1}^n x_k = 1$, we have  $1 \in H$ and $\lambda^N \in H$ for $N \in \mathbb{N}$. It is clear that $\frac{1}{\lambda^M}H \subset \frac{1}{\lambda^{M+1}}H$, $M \in \mathbb{N}$.
\end{remar}


One of the main results of this section is the following:

\begin{thm}\label{grouplike}
Let $\mu$ be an ergodic invariant measure on a stationary diagram $B$ defined by a distinguished eigenvalue $\lambda$ of the matrix $A = F^T$. Let $(x_1, \ldots, x_n)^T$ be the corresponding reduced vector and $H$  the additive subgroup of $\mathbb{R}$ generated by $\{x_1, \ldots , x_n\}$. Then the clopen values set $S(\mu)$ is group-like and $$
S(\mu) = \left(\bigcup_{N=0}^\infty \frac{1}{\lambda^N} H\right) \cap [0,1].
$$
\end{thm}

\noindent\textbf{Proof.} The proof is divided into two parts depending on the properties of $\lambda$. The first part deals with rational (hence integer) $\lambda$, and the second one contains the proof of the case of irrational (hence algebraic integer) $\lambda$.

\textbf{1.} Let $\lambda \in \mathbb{Q}$ and $x = (\frac{p_1}{q}, \ldots, \frac{p_n}{q})^T$ be the corresponding reduced probability vector. It follows from (\ref{formulaSmu}) that
$$
S(\mu) = \left\{\sum\limits_{i = 1}^{n} l^{(N)}_{i}\frac{p_{i}}{q \lambda^{N-1}}\; \vert\;  0 \leq l^{(N)}_{i} \leq h^{(N)}_{i},  \; N = 1,2,...\right\}.
$$
Hence, $S(\mu) \subset \{\frac{m}{q \lambda^N} :  N \in \mathbb{N},\; m =  0, 1,...,q \lambda^N\}$.
We need to prove the converse, that is, for every natural number $N$ and every integer $0 \leq m \leq q \lambda^N$, there exist $M \in \mathbb{N}$ and integers $l_i^{(M)} \in [0, h_{i}^{(M)}], i=1,...,n,$ such that
\begin{equation}\label{M}
\frac{m}{q \lambda^N} = \sum\limits_{i = 1}^{n} l^{(M)}_{i}\frac{p_{i}}{q \lambda^{M-1}}
\end{equation}
or, equivalently,
$$
m \lambda^R = \sum_{i = 1}^n l_i^{(M)}p_i,
$$
where $R = M - N - 1$.
We may assume that $0 < m < q \lambda^N$ because the cases $m = 0$ and $m = q \lambda^N$ are trivial. We note also that if an integer $M$ satisfying (\ref{M}) exists then  $M$ can be chosen arbitrary large, in particular, $M > N$.

Let $\alpha$ be the class of vertices corresponding to $\lambda$ and defining $\mu$.
If the measure $\mu$ is supported on a simple subdiagram of $B$ corresponding to a minimal component, then there is nothing to prove since this measure is invariant for a uniquely ergodic homeomorphism of a Cantor set.

Without loss of generality, we may assume that the non-zero value $\frac{p_1}{q}$ is assigned to a vertex from a class $\beta$ such that $\alpha \succ \beta$ (see (\ref{Frobenius Form})).

Since $\gcd(p_1, \ldots, p_n) = 1$, there exist integers  $d_1, \ldots, d_n$ such that $d_1 p_1 + \ldots + d_n p_n = 1$.
For the homogenous equation
$$
\sum_{i = 1}^n z_i p_i =0,
$$
there are $n-1$ independent parameters, say $(z_2,...,z_n)$,  amongst the solution of this equation. Then $z_1 = - \frac{1}{p_1} \sum_{i = 2}^n z_i p_i$. It is obvious that we can choose parameters $\{z_i\}_{i=2}^n$ such that all the numbers $\{z_i\}_{i=1}^n$ are integers.
From the above relations we  obtain that
$$
m \lambda^R = \sum_{i = 1}^n y_i p_i,
$$
where $y_i = m \lambda^R d_i + z_i$.
We need to show that there exist $z_2, \ldots, z_n$ such that $y_i \in \mathbb{N}$ and $0 \leq  y_i \leq h_{i}^{(M)}$ for $i = 1, \ldots, n$. We first note that
\begin{equation}\label{borders}
- m \lambda^R d_j \leq z_j \leq  h_{j}^{(N + R + 1)} - m \lambda^R d_j, \; j = 2,...,n.
\end{equation}
Since $y_1$ must be in the interval $[0,  h_{1}^{(M)}]$, the value $z_1$ must satisfy the inequalities
\begin{equation}\label{plus}
- m \lambda^R d_1 \leq z_1 \leq  h_{1}^{(N + R + 1)} - m \lambda^R d_1.
\end{equation}
On the other hand, it follows from (\ref{borders}) that
\begin{equation}\label{brd}
\frac{m \lambda^R}{p_1} \sum_{j = 2}^n d_j p_j - \sum_{j = 2}^n \frac{h_{j}^{(N + R+1)}p_j}{p_1}\leq z_1 \leq \frac{m \lambda^R}{p_1} \sum_{j = 2}^n d_j p_j.
\end{equation}
Since $\sum_{j = 2}^n d_j p_j = 1 - d_1 p_1$, we deduce from (\ref{brd})  that
\begin{equation}\label{minus}
\frac{m \lambda^R}{p_1} -  m \lambda^R d_1 - \frac{1}{p_1} \sum_{j = 2}^n h_{j}^{(N + R+1)}p_j \leq z_1 \leq \frac{m \lambda^R}{p_1} - m \lambda^R d_1.
\end{equation}
Thus, $z_1$ must satisfy both inequalities (\ref{plus}) and (\ref{minus}). We show that if $z_1$ satisfies (\ref{plus}) then $z_1$ also satisfies (\ref{minus}) when $R$ is sufficiently large.

To do this, we compare the left bounds of (\ref{minus}) and (\ref{plus}) and show that for sufficiently large $R$
\begin{equation}\label{p+-p-}
\frac{m \lambda^R}{p_1} - \frac{1}{p_1} \sum_{j = 2}^n h_{j}^{(N + R+1)}p_j = \frac{\lambda^{N + R}q}{p_1} \left(\frac{m}{q\lambda^N} - \sum_{j = 2}^n h_{j}^{(N + R+1)}\frac{p_j}{\lambda^{N + R}q}\right) < 0.
\end{equation}
Indeed, the vertex of the diagram corresponding to $\frac{p_1}{q}$  belongs to a final class of the vertices. Since
$$
\mu(X_B) = \sum_{j = 1}^n h^{(M)}_j \frac{p_j}{q \lambda^{M - 1}} = 1, \ \ M \in \mathbb{N},
$$
we have from asymptotics   (\ref{eq-asymp2})  that $h^{(M)}_1 \cdot \frac{p_1}{q \lambda^{M - 1}}\to 0$  as $M\to \infty$. Thus,
$$
\sum_{j = 2}^n h_{j}^{(N + R+1)}\frac{p_j}{\lambda^{N + R}q} \rightarrow 1
$$
as $R \rightarrow + \infty$.
Since $m < \lambda^{N}q$, the expression in parentheses in (\ref{p+-p-}) is negative for sufficiently large $R$ as desired. Moreover, the absolute value of expression (\ref{p+-p-}) tends to $+ \infty$ as $R \rightarrow + \infty$.

Similarly, comparing the right bounds of  (\ref{minus}) and (\ref{plus}), and using asymptotics $h_{1}^{(N + R+1)} = o(\lambda^R)$, we obtain that
\begin{equation}\label{left bounds}
\frac{m \lambda^R}{p_1} - h_{1}^{(N + R+1)} >0
\end{equation}
when $R$ is sufficiently large. Therefore, the interval defined by (\ref{plus}) lies in that  defined by (\ref{minus}) and its length tends to infinity as $R \rightarrow + \infty$.

To finish the proof, we need to show that $z_2,...,z_n$ can be chosen so that they satisfy simultaneously (\ref{borders}) and (\ref{plus}). We consider $z_1$ as a linear function of the parameters $z_2,...,z_n$. The integer parameters such that $z_1$ is integer and equations (\ref{borders}) hold form the domain of $z_1$. Relation (\ref{minus}) contains the range of $z_1$. The range of $z_1$ is a finite number of points.
The largest distance between two neighboring points is bounded and does not depend on $R$.
The length of subinterval (\ref{plus}) tends to infinity as $R$ becomes infinitely large. Hence, we can find allowable parameters $\{z_j\}_{j = 2}^n$ such that $z_1$ lies in the interval (\ref{plus}).
Thus, we set
$$
l_j^{(M)} = m \lambda^R d_j + z_j, \; j = 1,...,n
$$
where $M = R + N + 1$. This proves the theorem in the rational case.


\textbf{2.} Let $\lambda \in \mathbb{R} \setminus \mathbb{Q}$ and $x = (x_1,...,x_n)^T$ be the corresponding reduced vector.

To clarify the main idea of the proof,  we first consider an example. Let
$$
A =
\begin{pmatrix}
1 & 1 \\
1 & 2 \\
\end{pmatrix}.
$$
We have the eigenvalue $\lambda = \frac{3 + \sqrt{5}}{2}$ and eigenvector $x = (\frac{3-\sqrt{5}}{2}, \frac{\sqrt{5}-1}{2})^T$ for $A$. Hence, $x_1 = 3 - \lambda$, $x_2 = \lambda - 2$ and $\frac{1}{\lambda} = 3 - \lambda$. Then, by (\ref{formulaSmu}),
$$
S(\mu) = \{l_1^{(N)} (3 - \lambda)^N + l_2^{(N)} (\lambda - 2)(3 - \lambda)^{N-1} \; \vert \;0 \leq l^{(N)}_{i} \leq h^{(N)}_{i},\; i = 1,2; \; N = 1,2,...\}.
$$
The minimal polynomial for $\lambda$ is $f(t) = t^2 - 3t + 1$. It can be proved that $(3 - \lambda)^N = -f_{2N - 1} \lambda + f_{2N + 1}$, where $f_i$ is the $i$-th Fibonacci number. Hence $S(\mu)$ can be written in terms of polynomials of $\lambda$ of first degree:
$$
S(\mu) = \{l_1^{(N)} (-f_{2N-1}\lambda + f_{2N+1}) + l_2^{(N)}(f_{2N-2}\lambda - f_{2N}),N\geq 1\},
$$
where $0 \leq l^{(N)}_{i} \leq h^{(N)}_{i},\; i = 1,2$. Instead of polynomials, we can work with vectors formed by their  coefficients. Thus, we obtain a vector representation of any element from $S(\mu)$. Let $P_N = \{l_1^{(N)} (-f_{2N-1}\lambda + f_{2N+1}) + l_2^{(N)}(f_{2N-2}\lambda - f_{2N})\; \vert \;0 \leq l^{(N)}_{i} \leq h^{(N)}_{i},\; i = 1,2\}$. Then $P_N$ is a part of the lattice in $\mathbb{R}^2$ generated by vectors $u_N = (f_{2N+1},-f_{2N-1})^T$ and $v_N = (- f_{2N},f_{2N-2})^T$ which includes all points with coordinates $\{(i,j)\; \vert \; 0 \leq i \leq h^{(N)}_{1}, 0 \leq j \leq h^{(N)}_{2}\}$ in the basis $\{u_N,v_N\}$. We see that $P_{N+1} \supset P_N$ and $S(\mu) = \bigcup_{i=1}^{\infty} P_N$. It can be proved that $u_N, v_N$ both tend to the same line $a$ in $\mathbb{R}^2$ generated by the vector $(-\lambda,1)^T$ as $N\to \infty$. We also show that the norms of these vectors tend to infinity. This suffices to conclude that the points in $\mathbb{R}^2$ that represent $S(\mu)$ ``uniformly'' fill the gap between lines $a$ and $a + (1,0)^T$. This means that $S(\mu)$ is group-like.
\medskip

Now we consider the general case.
It suffices to prove that for any $u, v \in S(\mu)$ with $u + v \leq 1$ we have $u + v \in S(\mu)$. Indeed, let $u, v \in S(\mu)$ with $u < v$. Then $v - u = 1 - ((1 - v) + u) \in S(\mu)$ (Akin used similar arguments in \cite{Akin2}). Then  $S(\mu)$ is group-like by Lemma \ref{AkinKritGrouplike}. It follows from  (\ref{formulaSmu}) that it suffices to prove that any number $s = \sum_{i=1}^n l_i \frac{x_i}{\lambda^{N-1}}$ from $[0,1]$ such that $l_i \geq 0$ belongs to $S(\mu)$.

We will use a vector representation of algebraic numbers as we did in the above example.  Since $\lambda$ is a root of the  characteristic polynomial  $\Delta(x)$ of $A$, $\lambda$ is an algebraic integer number. Suppose the degree of $\lambda$ is $k$. Denote by
$$
f(t) = t^k + m_{k-1}t^{k-1}+ \ldots+m_1 t + m_0
$$
the minimal polynomial of $\lambda$ over $\mathbb{Q}$. Let  $\mathbb{Q(\lambda)}$ denote the least field that contains both $\mathbb{Q}$ and $\lambda$.
If  $\mathbb{Q[\lambda]}$ stands for the least ring that contains both $\mathbb{Q}$ and $\lambda$, then $\mathbb{Q[\lambda]} = \mathbb{Q(\lambda)}$.
Any element  $y\in \mathbb{Q(\lambda)}$ can be uniquely represented as $a_0 + a_1 \lambda + \ldots + a_{k-1} \lambda^{k-1} = y(\lambda) $ where $a_i \in \mathbb{Q}$.  Hence, there is a one-to-one correspondence between polynomials in $\lambda$ and vectors formed by  their coefficients:
\begin{equation}\label{Correspondence}
a_0 + a_1 \lambda + \ldots + a_{k-1} \lambda^{k-1} \leftrightarrow
(a_0, a_1,...,a_{k-1})^T.
\end{equation}
Since $S(\mu) \subset \mathbb{Q}(\lambda)$, every  element of $S(\mu)$ can be also considered as a vector in the space $\mathbb{Q}^k$. Thus, if we need to emphasize that a number $y \in S(\mu)$ is considered as a vector from $ \mathbb{Q}^k$   (or a  polynomial from  $\mathbb{Q}(\lambda)$), we will use the notation $\textbf y$. This convention will be used throughout the paper.

It follows from (\ref{Correspondence})  that the polynomials $\{1, \lambda,...,\lambda^{k-1}\}$ correspond to the vectors $\{\textbf{e}_1, \textbf{e}_2,...,\textbf{e}_k\}$,  the standard basis in $\mathbb{R}^k$ (or $\mathbb{Q}^k$). Denote by $\langle \textbf{u}, \textbf{v}\rangle$ the scalar product of vectors $\textbf{u},\textbf{v} \in \mathbb{R}^k$.
Let $\textbf{n}$ denote the vector $(1, \lambda, ..., \lambda^{k-1})^T \in \mathbb{R}^k$. Then, for any $y\in \mathbb Q^k$, the corresponding polynomial $y(\lambda)$ can be written as $y(\lambda) = \langle \textbf{y}, \textbf{n}\rangle \in \mathbb{R}$.

Take the entries $x_1,...,x_n$ of the reduced vector, and for each $x_i$ find its representation as a vector from $\mathbb{Q}^k \subset \mathbb{R}^k$:
$$
\textbf{x}_i = (a_0^{(i)}, a_1^{(i)},...,a_{k-1}^{(i)})^T.
$$
Denote by
$$
A_0 =
\begin{pmatrix}
a_0^{(1)} & \ldots & a_0^{(n)}\\
a_1^{(1)} & \ldots & a_1^{(n)}\\
\vdots & \ldots & \vdots \\
a_{k-1}^{(1)} & \ldots & a_{k-1}^{(n)}
\end{pmatrix}
$$
the matrix formed by the vectors $\textbf{x}_1,...,\textbf{x}_n$.
In other words, $A_0$ represents the transposed eigenvector $(x_1,\ldots, x_n)$.

Let $p(\lambda)$ be the polynomial in $\mathbb Q(\lambda)$ such that $\lambda^{-1}= p(\lambda)$. The map $y \mapsto p(\lambda)y$ in $\mathbb Q(\lambda)$ determines a linear transformation in the vector space $\mathbb{Q}^k$. Find the $k\times k$ matrix $D$ which corresponds to  this transformation.  It is obvious that $D({\textbf e_i}) = {\textbf e_{i-1}}, i=2,...,k,$ and $D({\textbf e_1})$ is the vector corresponding to $p(\lambda)$. Thus, the matrix $D$ and the inverse matrix $C = D^{-1}$ have the form
$$
D =
\begin{pmatrix}
-\frac{m_1}{m_0} & 1 & \ldots & 0\\
\vdots & \vdots & \ddots & \vdots \\
-\frac{m_{k-1}}{m_0} & 0 & \ldots & 1\\
-\frac{1}{m_0} & 0 & \ldots & 0
\end{pmatrix}, \ \ \
C =
\begin{pmatrix}
0 & \ldots & 0 & - m_0\\
1 & \ldots & 0 & - m_1\\
\vdots & \ddots & \vdots & \vdots \\
0 & \ldots & 1 & - m_{k-1}\\
\end{pmatrix},
$$
where  $m_0, m_1,...,m_{k-1}$ are the coefficients of the minimal polynomial $f(t)$.
Then the equation $Ax = \lambda x$ can be written in new terms as $C A_0 = A_0 A^T$ or $A_0 = D A_0 A^T$.

\medskip

Now we use the found matrix $D$ to represent the clopen values set in more convenient form:
$$
S(\mu) = \left\{D^{N-1} \left(\sum_{i = 1}^{n} k^{(N)}_{i} {\textbf x}_i\right) \ |\  0  \leq k^{(N)}_{i} \leq h^{(N)}_{i}; \; N = 1,2, ...\right\} \subset \mathbb{Q}^k.
$$
Denote by $\pi$ the hyperplane in $\mathbb{R}^k$ which is specified by the relation $\pi = \{\textbf y : \langle \textbf{y}, \textbf{n}\rangle = 0\}$. Let $\pi_1$ be the hyperplane which is obtained by the relation $\pi_1 = \{\textbf y : \langle \textbf{y}, \textbf{n}\rangle = 1\}$ (hence $\pi_1 = \pi + \textbf{e}_1$). Then all points of $S(\mu)$ lie in the stripe between $\pi$ and $\pi_1$.

Consider the family of sets
$$
P_N = \left\{\sum_{i=1}^n l_i^{(N)} D^{N-1} (\textbf{x}_i) \ | \  0 \leq l_i^{(N)} \leq h_i^{(N)} \right\},\ \ N \in \mathbb{N}.
$$
Since
$$
\frac{x_i}{\lambda^{N}} = \sum_{i=1}^n a_{ij}\frac{x_j}{\lambda^{N+1}},
$$
we see that $P_{N+1} \supset P_N$ for $N \in \mathbb{N}$. Let $P = \bigcup_{N \in \mathbb{N}} P_N$. Clearly, $S(\mu) = P$.

We show that the vectors from $P$ can be chosen arbitrary close to the hyperplane $\pi$ and have arbitrary big length. Hence they ``fill'' the gap between hyperplanes $\pi$ and $\pi_1$.
Thus we show that for any two vectors $\textbf{u},\textbf{v} \in P$ if $\textbf{u} + \textbf{v}$ lies in the stripe between hyperplanes $\pi$ and $\pi_1$ then there exists $N \in \mathbb{N}$ such that $\textbf{u}+\textbf{v} \in P_N \subset P$.

Recall that the non-zero entries $x_1,...,x_n$ are related to the vertices that are accessible from the distinguished class $\alpha$ defining the measure $\mu$. Suppose that the vertices $m+1, \ldots, n$ belong to $\alpha$. We will need the following lemma.

\begin{lemm}\label{finalbasis}
(1) $\det (C - \lambda I) = 0$ if and only if $f(\lambda) = 0$. All eigenvalues of $C$ (and hence $D$) are distinct.

(2) The rank of the matrix $A_0$ is $k$. Moreover, a set of $k$ linearly independent columns of $A_0$  can be chosen amongst the vectors $\{\textbf{x}_{m+1},...,\textbf{x}_n\}$.
\end{lemm}

\noindent\textbf{Proof.}
In order to prove (1), it suffices to notice that $\det(C - tI) = (-1)^k f(t)$.

Because the eigenvector $(x_1,...,x_n)^T$ is probability, we obtain
that for any $M\in \mathbb N$
\begin{equation}\label{lambda^M}
\lambda^M = \sum_{i=1}^n \sum_{j=1}^n (A^M)_{ij}x_j = \sum_{j=1}^n x_j \sum_{i=1}^n(A^M)_{ij},
\end{equation}
and therefore every basis vector $\textbf{e}_i$ can be represented as a linear combination of vectors $\textbf{x}_1,...,\textbf{x}_n$ with non-negative coefficients.
Let $S$ be the subspace $\mathbb{Q}^k$ generated by vectors $\{\textbf{x}_{1},...,\textbf{x}_n\}$. First we show that $\dim(S) \geq k$. Indeed, relation (\ref{Correspondence}) shows that the numbers $1, \lambda,...,  \lambda^{k-1}$ correspond to the standard basis in $\mathbb{Q}^k$. Then we deduce from (\ref{lambda^M}) that this basis belongs to $S$ and, therefore, $\rm{rank}(A_0) = k$.

Next, we use the fact that $\lambda$ is also the Perron-Frobenius eigenvalue of the submatrix $A_\alpha = F_\alpha^T$. Applying (\ref{lambda^M}) to $A_\alpha$, we obtain that
\begin{equation}\label{rankA_0}
\sum_{i=m+1}^n \sum_{j=m+1}^n (A^r)_{ij}x_j = \lambda^r \sum_{i=m+1}^n x_i
\end{equation}
for $r \in \mathbb{N}$. We denote  $q(\lambda) = \sum_{i=m+1}^n x_i$. Then $q(\lambda) \in  \mathbb{Q}[\lambda]$. It follows from (\ref{rankA_0}) that, for any polynomial $h(\lambda) \in  \mathbb{Q}[\lambda]$, the product $h(\lambda)q(\lambda)$ is a linear combination of vectors $\{\textbf{x}_{m+1},...,\textbf{x}_n\}$  with rational coefficients. Then,  choose $h_r\in \mathbb{Q}[\lambda]$ such that
 $\lambda^r = h_r(\lambda)q(\lambda)$. Therefore, each $\lambda^r, \ r=0,...,k-1$, can be represented as a linear combination of $\{\textbf{x}_{m+1},...,\textbf{x}_n\}$  with rational coefficients. It follows that the set $\{\textbf{x}_{m+1},...,\textbf{x}_n\}$  contains $k$ linearly independent vectors.  \hfill$\blacksquare$
\medskip

We continue the proof of the theorem. We note that $D^T \textbf{n} = \frac{1}{\lambda}\, \textbf{n}$. It is clear that $D(\pi) = C(\pi) = \pi$.

Since $(-1)^{k} f(x)$ is the characteristic polynomial for $C$, the matrix $C$ has $k$ distinct eigenvalues in $\mathbb{C}$. In particular, $\lambda$ is an eigenvalue for $C$, and for every other eigenvalue $\nu$ we have $|\nu| < \lambda$. It is not hard to find the eigenvector $\textbf{y}_1$ of $C$ corresponding to  $\lambda$:
$$
\textbf{y}_1 =
\begin{pmatrix}
1\\
\frac{1}{\lambda} + \frac{m_1}{m_0}\\
\vdots\\
\left(\frac{1}{\lambda}\right)^{k-1} + \frac{m_1}{m_0}\left(\frac{1}{\lambda}\right)^{k-2} + ... + \frac{m_{k - 1}}{m_0}
\end{pmatrix}.
$$
This is the only eigenvector of $C$ (and hence of $D$) which is not orthogonal to $\textbf{n}$. Indeed, suppose $C\textbf{z} = \nu \textbf{z}$. Then
$$
\nu\langle \textbf{z}, \textbf{n}\rangle = \langle C \textbf{z}, \textbf{n}\rangle = \langle \textbf{z},C^T \textbf{n} \rangle = \lambda\langle \textbf{z}, \textbf{n}\rangle.
$$
We see that either $\nu = \lambda$ or $\langle \textbf{z}, \textbf{n}\rangle = 0$. All eigenvalues of $C$  are different, so that we obtain that $\langle \textbf{z}, \textbf{n}\rangle = 0$. This means that all eigenvectors of $C$ (or $D$) different from $\textbf{y}_1$ belong to $\pi$.

Now we are ready to prove that the vectors of $P$ can be chosen arbitrary close to the hyperplane $\pi$. Consider two cases when the eigenvalues of $C$ are real or complex.

(I) Suppose all eigenvalues of $C$ are real. Then $C$ is diagonalizable in $\mathbb{R}^k$. Let $\textbf{y}_1,...,\textbf{y}_k$ be the eigenvectors of $C$ which correspond to the eigenvalues $\lambda_1 = \lambda,...,\lambda_k$. For $\textbf{y} \in \mathbb{R}^k$, there exist unique real numbers $\alpha_i$ such that $\textbf{y} = \sum_{i=1}^k \alpha_i \textbf{y}_i$. Then $C^N \textbf{y} = \alpha_1 \lambda^N \textbf{y}_1 + \textbf{z}_N$ where $\textbf{z}_N = \sum_{i=2}^k \alpha_i \lambda_i^N \textbf{y}_i$ belongs to  $\pi,\ N \in \mathbb{N}$. We prove that the angle between the lines generated by $C^N \textbf{y}$ and $\textbf{y}_1$ can be made arbitrary small when $N$ tends to infinity. Indeed, it is obvious that
$$
\frac{\parallel\textbf{z}_N\parallel}{\parallel \alpha_1 \lambda^N \textbf{y}_1\parallel} =  \frac{\sqrt{\sum_{2 \leq i,j \leq k}\alpha_i \alpha_j \langle \textbf{y}_i, \textbf{y}_j \rangle \lambda_i^N \lambda_j^N}}{\parallel \alpha_1 \lambda^N \textbf{y}_1\parallel} \rightarrow 0
$$
as  $N \rightarrow \infty$ because  $|\lambda_i| < \lambda$ for $i = 2,...,n$.

(II) Suppose some of the eigenvalues of $C$ are complex. Then they form the pairs of complex conjugate numbers; say, $\lambda_2 = \alpha + i \beta$, and $\lambda_3 = \alpha - i \beta$ where $\alpha, \beta  \in \mathbb{R}$. Let $\textbf{u} + i \textbf{v}$ be the eigenvector of $C$ corresponding to $\lambda_2$ where $\textbf{u}, \textbf{v} \in \mathbb{R}^k$. Then the subspace of $\mathbb{R}^k$ generated by $\textbf{u},\textbf{v}$ is an invariant space for $C$ and $C\textbf{u} = \alpha \textbf{u} - \beta \textbf{v}$, $C\textbf{v} = \alpha \textbf{v} + \beta \textbf{u}$. Since $|\lambda_2| < \lambda$, we have $|\alpha| < \lambda$ and $|\beta| < \lambda$.

We represent $\textbf{y}$ as a linear combination of real eigenvectors and real components of complex eigenvectors of $C$. The proof now is analogous to that in case (I).

Thus, while the iterations of $C$ drive any ray which is not in $\pi$ to the limit ray generated by $\textbf{y}_1$, the iterations of $D = C^{-1}$ do the opposite thing. Arguing as above, we can prove that the angle between the line generated by $D^N \textbf{y}$ and $\pi$ can be made arbitrary small when $N$ tends to infinity.

Applying the iterations of $D$, we prove that the length of vectors from $P$ can be made arbitrary long. Recall that the vertices $m+1, \ldots, n$ belong to the distinguished class $\alpha$ corresponding to $\mu$. By Lemma \ref{finalbasis}, the vectors $\{x_i\}_{i=m+1}^n$ contain a basis of $\mathbb{R}^k$.
We have $\langle D^N \textbf{x}_i, \textbf{n}\rangle \rightarrow 0$ as $N \rightarrow \infty$. Moreover, $\cos\angle(D^N \textbf{x}_i, \textbf{n}) \rightarrow 0$ as $N \rightarrow \infty$. Consider $\langle h^{(N+1)}_i D^N \textbf{x}_i, \textbf{n}\rangle$. For $i = m+1,...,n$, it follows from the asymptotics of $(A^N)_{ij}$ that $h^{(N+1)}_i \frac{x_i}{\lambda^N} \sim c_i x_i$ as $N \rightarrow \infty$ where $c_i$ is a constant. Hence $\langle h^{(N+1)}_i D^N \textbf{x}_i, \textbf{n}\rangle \sim c_i x_i$ and finally $\| h^{(N+1)}_i D^N \textbf{x}_i \| \rightarrow \infty$ as $N \rightarrow \infty$ because $\cos\angle(D^N \textbf{x}_i, \textbf{n}) \rightarrow 0$.

Recall that we need to prove that any number $s = \sum_{i=1}^n l_i \frac{x_i}{\lambda^{N-1}}$ from $[0,1]$ such that $l_i \geq 0$ lies in $S(\mu)$. In the vector interpretation,  this means that $\textbf{s} = \sum_{i=1}^n l_i D^{N-1} \textbf{x}_i$. We will use in our proof the fact that the measure $\mu$ is  supported on the set $X_\alpha$ of all infinite paths that eventually go through  vertices of the class $\alpha$. We consider two cases.

(i)  Let $0 \leq l_i \leq h_i^{(N)}$ for $i \not \in \alpha$ and $0 < s =  \langle \textbf{s}, \textbf{n} \rangle < 1$. We show that $s \in S(\mu)$. Indeed, the part $\textbf s' = \sum_{i\not \in \alpha} l_i D^{N-1} \textbf{x}_i$ belongs to $P_N$ because its coefficients $l_i$ lie in the needed range. Clearly,  $\textbf{s}' \in P_M$ for $M > N$. The vector  $\sum_{i \in \alpha} l_i D^{N-1} \textbf{x}_i$ lies in the integer lattice generated by $D^M \textbf{x}_{m+1},...,D^M \textbf{x}_m$ for $M \geq N$. Since $\| h_i^{(M)}D^{M-1}\textbf{x}_i \| \rightarrow \infty$ as $M\rightarrow  \infty$, $i = m+1,...,n$,  the allowable linear combinations of $D^M \textbf{x}_{m+1},...,D^M \textbf{x}_m$ eventually fill the stripe between $\pi$ and $\pi_1$. Thus, the point $\textbf{s}$ will be ``covered'' by an  allowable combination.

(ii) Let $\textbf{s} = \sum_{i=1}^n l_i^{(N)} D^{N-1} \textbf{x}_i$, where $l_i^{(N)} \in \mathbb{N}$ and $0 < \langle \textbf{s}, \textbf{n} \rangle < 1$. We show that this case can be reduced to the previous one. Suppose there exist $l_i > h_i^{(N)}$ for some $i \in \{1,...,m\}$. By Lemma \ref{finalbasis}, there exist  $\{p_{ij}\}_{j = m + 1}^n \in \mathbb{Z}$ and $q_i \in \mathbb{N}$ such that
\begin{equation}\label{combi}
q_i \textbf{x}_i = \sum_{j = m+1}^n p_{ij} \textbf{x}_j
\end{equation}
for $i  = 1,...,m$. We can find $M > N$ such that $h_i^{(M)} > q_i$ for $i  = 1,...,m$. We have $\textbf{s} = \sum_{i=1}^n l_i^{(M)} D^{M-1} \textbf{x}_i$ where $l_i^{(M)} = \sum_{j=1}^n l_j^{(N)} (A^{(M - N)})_{ji}$. Then if $l_i^{(M)} > h_i^{(M)}$ for some $i \in \{1,...,m\}$ we write down $l_i^{(M)} = t_i^{(M)} q_i + r_i^{(M)}$ where $t_i^{(M)} \in \mathbb{N}$ and $0 \leq r_i^{(M)} < q_i$. By (\ref{combi}), the vector $t_i^{(M)} q_i \textbf{x}_i$ can be expressed as the integer combination of $\{\textbf{x}_j\}_{j=m+1}^n$. Hence
$$
\textbf{s} = \sum_{i=1}^m r_i^{(M)} D^{M-1} \textbf{x}_i + \sum_{i=m+1}^n \left(l_i^{(M)} + \sum_{j=1}^m p_{ji} t_j^{(M)}\right)D^{M-1} \textbf{x}_i.
$$
Since $r_j^{(M)} < q_j \leq h_j^{(M)}$ for $j = 1,...,m$, it suffices to show that the coefficients $\left(l_i^{(M)} + \sum_{j=1}^m p_{ji} t_j^{(M)}\right)$, $i = m+1,...,n$, can be made positive for $M$ large enough. It follows from the above relations that  $l_i^{(M)} \sim \lambda^M$ as $M \rightarrow \infty$ for $i \in \alpha$. On the other hand, $t_j^{(M)} < l_j^{(M)}$ and  $l_j^{(M)} \sim \bar{o}(\lambda^M)$ as $M \rightarrow \infty$ for $j = 1,...,m$. Hence the needed coefficients can be made positive. \hfill$\blacksquare$

\medskip

\textbf{3.2. Good measures.} Now we consider the conditions under which an ergodic invariant measure on a stationary Bratteli diagram is good.

\begin{lemm}\label{gdfin}
Let $\mu$ be an ergodic $\mathcal R$-invariant measure on a stationary Bratteli diagram $B$ and let $A$ be the matrix transposed to the incidence matrix of $B$. Denote by $\alpha$  the distinguished class of vertices that defines  $\mu$. Then $\mu$ is good if and only if all the clopen cylinder sets that end in the vertices of the class $\alpha$ are good.
\end{lemm}

\noindent \textbf{Proof.}
The "only if" part of this result is obvious.

To prove the "if" part, we consider any clopen sets $U, V \subset X_\alpha$  with $\mu(U) < \mu(V)$. Since cylinder sets form a partition basis for $X_B$, we may assume that $V$ is a cylinder set (see Theorem \ref{goodmeasure}). We must find a clopen subset $W \subset V$ with $\mu(U) = \mu(W)$.
By definition of $X_\alpha$, the finite path corresponding to  $V$ ends in a  vertex of a class $\beta$ that is accessible from $\alpha$, i.e.  $\alpha \succeq  \beta$  (otherwise $V$ would have zero measure). If $\beta = \alpha$, then there is nothing to prove.
Suppose now that $\alpha \succ \beta$. Denote by $N$ the length of the cylinder set $V$. Then $V$ is a disjoint union of cylinder subsets of length $N + 1$. Their end vertices are either in the class $\alpha$ or in the classes that are accessible from  $\alpha$. Take the latter cylinder sets and represent each of them  as a disjoint union of cylinder sets of length $N + 2$.
We continue these partitions infinitely many times. Enumerate the cylinder subsets of $V$ that end in vertices of $\alpha$ by $V_k,\ k\in \mathbb N$. Then
$$V \supset \bigsqcup\limits_{k=1}^{\infty} V_k \quad \mbox{and}\quad
\mu(V) = \sum_{k=1}^{\infty} \mu(V_k)
$$
by asymptotics (\ref{eq-asymp}), (\ref{eq-asymp2}). Since $\mu(U) < \mu(V)$, we can find  $m$ such that
$$
\sum_{k=1}^{m}\mu(V_k) \leq \mu(U) < \sum_{k=1}^{m+1}\mu(V_k).
$$
If the equality $\sum_{k=1}^{m}\mu(V_k)= \mu(U)$ holds, then we define  $W = \bigsqcup\limits_{k=1}^{m} V_k$. Otherwise consider the number $\mu(U) - \sum_{k=1}^{m}\mu(V_k)$. This number is contained in $S(\mu)$, since $S(\mu)$ is group-like by Theorem \ref{grouplike}. The set $V_{m+1}$ is good because it ends in a vertex of $\alpha$. Hence, it contains a clopen subset $U_1$  such that $\mu(U_1) = \mu(U) - \sum_{k=1}^{m}\mu(V_k)$. Then we set $W = \bigsqcup\limits_{k=1}^{m} V_k \bigsqcup U_1$. \;\hfill$\blacksquare$

\begin{thm}\label{KritGood}
Let $\mu$ be an ergodic $\mathcal R$-invariant measure on a stationary diagram $B$ defined by a distinguished eigenvalue $\lambda$ of the matrix $A = F^T$. Denote by $x = (x_1,...,x_n)^T$ the corresponding reduced vector. Let the vertices $m+1, \ldots, n$ belong to the distinguished class $\alpha$ corresponding to $\mu$. Then $\mu$ is good if and only if   there exists $R \in \mathbb{N}$ such that $\lambda^R x_1,...,\lambda^R x_m$ belong to the additive group generated by $\{x_j\}_{j=m+1}^n$.

If the clopen values set of $\mu$ is rational and $(\frac{p_1}{q}, \ldots, \frac{p_n}{q})^T$ is the corresponding reduced vector, then $\mu$ is good if and only if $\gcd(p_{m+1}, ...,p_n) | \; \lambda^R$ for some $R \in \mathbb N$.
\end{thm}

\noindent \textbf{Proof.}
If $m = 0$, then the Bratteli diagram $B$ is simple and the measure $\mu$ is good by Theorem \ref{goodmeasure}. Suppose $m > 0$. Consider the $(n-m) \times (n-m)$ block $A_\alpha$ of the matrix $A$ whose entries count the edges between  vertices of the class $\alpha$. Set $\widetilde{x} = \sum_{k = m + 1}^{n} x_k$. Then $\left(\dfrac{x_{m+1}}{\widetilde{x}}, \ldots, \dfrac{x_n}{\widetilde{x}}\right)^T$ is the probability eigenvector for $A_\alpha$ corresponding to the eigenvalue $\lambda$. Let $B_\alpha$ be the stationary subdiagram of $B$ consisting of vertices from the class $\alpha$ and edges connecting them. Then $A_\alpha$ is the matrix transpose to the incidence matrix of the subdiagram $B_\alpha$. Moreover, we can assume, without loss of generality, that  $B_\alpha$ is a simple subdiagram (see \cite{S.B.}).
Let $\widetilde{\mu}$ be the (unique) ergodic $\mathcal R$-invariant measure on this diagram.
If  $Y_\alpha$ denotes the path space of the Bratteli diagram $B_\alpha$, then $\mu(Y_\alpha) = \widetilde{x}$ because $Y_\alpha$ is a complete section in $X_\alpha$ for the tail equivalence relation. Then $\widetilde{\mu}$ can be regarded as a relative measure on the clopen subset $Y_\alpha$ of $X_\alpha$. By Theorem \ref{goodmeasure}, the measure $\widetilde{\mu}$ is good.

Denote by $H\left(\dfrac{x_{m+1}}{\widetilde{x}},...,\dfrac{x_{n}}{\widetilde{x}}\right)$ the group generated by $\dfrac{x_{m+1}}{\widetilde{x}},...,\dfrac{x_{n}}{\widetilde{x}}$. By Theorem~\ref{grouplike},
\begin{equation}\label{Smureduced}
S(\widetilde{\mu}) = \left\{\frac{r}{\lambda^N} : N \in \mathbb{N},\; 0 \leq r \leq \lambda^N, r \in H\left(\frac{x_{m+1}}{\widetilde{x}}, \ldots, \frac{x_n}{\widetilde{x}}\right)\right\}.
\end{equation}

Suppose that for any $x_i$, $i = 1,...,m$ there exists $R \in \mathbb{N}$ such that $x_i$ belongs to the additive group generated by $\left\{\frac{x_j}{\lambda^R}\right\}_{j=m+1}^n$. Then we can find one common $R$ for them.
By Lemma \ref{gdfin}, it suffices to prove that any cylinder set that ends in the class $\alpha$ is good. Let $V$ be a cylinder set whose terminal vertex $v$ is in the class $\alpha$. For any clopen set $U$ with $\mu(U) < \mu(V)$ we must find a clopen subset $W \subset V$ such that $\mu(W) = \mu(U)$. Consider a cylinder set $\widetilde{V}$ that ends at the same vertex $v$, but
passes only through the vertices of the class $\alpha$ (i.e. $\widetilde{V}$ is tail equivalent to $V$).
Then $\widetilde{V} \subset Y_\alpha$ and $\mu(\widetilde{V}) = \mu(V)$. To prove the claim, it suffices to find a clopen set $\widetilde{W} \subset \widetilde{V}$ such that $\mu(\widetilde{W}) = \mu(U)$. Indeed, if such a set $\widetilde{W}$ exists then we can find a tail equivalent clopen subset  $W \subset V$ such that $\mu(\widetilde{W}) = \mu(W)= \mu(U)$.
By Theorem \ref{grouplike}, $\mu(V) = \frac{x_j}{\lambda^N}$, for some $N \in \mathbb N$ and $j \in \{m+1, \ldots, n\}$. Since $U$ is a clopen set, we have $\mu(U) = \frac{k}{\lambda^M}$ where $k \in H$. We can take $M \geq N$. For any subset $\widetilde{W} \subset Y_\alpha$, we see that $\mu(\widetilde{W}) = \widetilde{x}\cdot \widetilde{\mu}(\widetilde{W})$. By (\ref{Smureduced}), $\widetilde{\mu}(\widetilde{W}) = \frac{k_1}{\lambda^S}$ for some $S \in \mathbb N$, and $0 \leq k_1 \leq \lambda^S$, $k_1 \in H\left(\frac{x_{m+1}}{\widetilde{x}}, \ldots, \frac{x_n}{\widetilde{x}}\right)$. Thus, we need to find $S, k_1$  such that
\begin{equation}\label{relation}
\mu(U) = \frac{k}{\lambda^M} = \frac{k_1}{\lambda^S} \cdot \widetilde{x} = \mu(\widetilde{W}).
\end{equation}
Let $k = \sum_{i=1}^n d_i x_i$ and $k_1 = \sum_{i=m+1}^n c_i \dfrac{x_i}{\widetilde{x}}$.
Then it follows from~(\ref{relation})
that
$$
\sum_{i=1}^n d_i x_i = \frac{1}{\lambda^{S-M}}\sum_{j=m+1}^n c_j x_j .
$$
Since $x_i \in H\left(\frac{x_{m+1}}{\lambda^{S-M}},...,\frac{x_{n}}{\lambda^{S-M}}\right)$ for $S-M \geq R$, there exist integers  $c_{m+1},...,c_n$ satisfying the above equation. Because $\widetilde{\mu}$ is good, we can find a clopen subset $\widetilde{W} \subset \widetilde{V}$ with $\widetilde{\mu}(\widetilde{W}) = \frac{k_1}{\lambda^S}$, i.e.  $\mu(\widetilde{W}) = \mu(U)$.

Conversely, suppose that $\mu$ is a good measure. We can repeat the proof of the "if" part backwards to obtain the needed result.

In the rational case, $x_i = \frac{p_i}{q}$, $H = \frac{1}{q}\mathbb{Z}$. Then
$H(x_{m+1},...,x_n) = \frac{a}{q}\mathbb{Z}$, where $a = \gcd(p_{m+1}, ...,p_n)$. Hence, $\mu$ is good if and only if $\gcd(p_{m+1}, ...,p_n) | \; \lambda^R$ for some $R \in \mathbb{N}$.
\hfill$\blacksquare$

\begin{corol}\label{gcd1}
Let $\mu$ be the rational measure (i.e. $S(\mu) \subset \mathbb Q$) on a stationary diagram $B$ defined by a  distinguished eigenvalue $\lambda$ of the matrix $A = F^T$. Denote by $(\frac{p_1}{q}, \ldots, \frac{p_n}{q})^T$ the corresponding reduced vector. Let the vertices $m+1, \ldots, n$ belong to the distinguished class $\alpha$ corresponding to $\mu$. If $\gcd(p_{m+1}, ...,p_n) = 1$, then $\mu$ is good.
\end{corol}

From Theorems \ref{goodmeasure} and \ref{grouplike} we obtain the following
\begin{corol}
For an ergodic invariant measure $\mu$ on a stationary Bratteli diagram the following are equivalent:

(i) $\mu$ is a good measure.

(ii) $\mu$ is refinable.

(iii) $\mu$ is weakly refinable.

\end{corol}
\medskip


\section{Homeomorphic Measures on Stationary Diagrams}\label{section4}

Let $\mathcal D$ be the set of all (non-simple) stationary Bratteli diagrams. Denote by $\mathcal S$ the set of all Borel probability measures on diagrams from $\mathcal D$ which are  ergodic and invariant with respect to the tail equivalence relation. Recall that, in other words, a measure $\mu\in \mathcal S$ if and only if there exist an aperiodic substitution dynamical system $(Y,\varphi)$ and an ergodic $\varphi$-invariant measure $\nu$ such that $\mu$ is homeomorphic to $\nu$ \cite{BKM}. Clearly, $\mathcal S$ is a countable set.

Our goal is to show that for every good measure $\mu$ from $\mathcal S$ there are countably many measures $\mu_i\in \mathcal S$ on stationary Bratteli diagrams $B_i$ such that $\mu_i$ is homeomorphic to $\mu, \ i\in \mathbb N$. Moreover, the stationary Bratteli diagrams $B_i$ can be chosen essentially different:  the corresponding tail equivalence relations $\mathcal R_i$ are pairwise non-orbit equivalent.

\begin{thm}\label{goodhomeo}
Let $\mu$ be a good ergodic $\mathcal{R}$-invariant probability measure on a stationary (non-simple) Bratteli diagram $B$. Then there exist stationary Bratteli diagrams  $\{B_i\}_{i=0}^\infty$ and  good ergodic $\mathcal{R}_i$-invariant probability measures $\mu_i$ on $B_i$ such that each measure $\mu_i$ is homeomorphic to $\mu$ and the dynamical systems $(B_i, \mathcal R_i)$, $(B_j, \mathcal R_j)$ are topologically orbit equivalent if and only if $i = j$. Moreover, the diagram $B_i$ has exactly $i$ minimal components for the tail equivalence relation $\mathcal R_i, i\in \mathbb N$.
\end{thm}

\noindent \textbf{Proof.} We divide the proof into two cases: (1) the set  $S(\mu)$ has only rational values; (2) there are irrational values in $S(\mu)$.

\textbf{1}. Let $S = S(\mu)\subset \mathbb Q$. Then, as proved in Theorem \ref{grouplike}, there exist natural numbers  $\lambda$ and $q$, greater than one, such that $S = \{\frac{m}{q \lambda^N}\; |\; m, N \in \mathbb{N},\; 0 \leq m \leq q \lambda^N\}$.

We first construct a simple Bratteli diagram $B_0$ and an ergodic probability invariant measure $\mu_0$ such that $S(\mu_0) = S$. H. Yuasa~\cite{Yuasa} used similar arguments in the study of  orbit equivalence of substitution systems arising from primitive substitutions whose composition matrices have rational Perron-Frobenius eigenvalues.
We take the probability vector $x = (\frac{1}{q},...,\frac{1}{q})^T$ and the $q\times q$ matrix
$$
A_0 =
\begin{pmatrix}
\lambda - 1 & 1 &  \ldots & 0 & 0\\
0 & \lambda - 1 &  \ldots & 0 & 0\\
\vdots & \vdots & \ddots & \vdots& \vdots \\
0 & 0 &  \ldots & \lambda - 1 & 1\\
1 & 0 &  \ldots & 0 & \lambda - 1
\end{pmatrix}.
$$
Then $A_0x = \lambda x$. Clearly, the stationary Bratteli diagram $B_0$ defined by the transpose to $A_0$ is simple. The unique ergodic probability measure $\mu_0$ is good and $S(\mu_0) = S$. Hence, $\mu$ and $\mu_0$ are homeomorphic by Corollary \ref{homeomorphic_measures}.

Next, fix $i\ge 1$ and construct Bratteli diagrams $B_i$ and measures $\mu_i$ as follows. Set $\lambda_i = \lambda^{i+1}$ and define the probability vector $x_i =(\frac{1}{q\lambda_i},..., \frac{1}{q\lambda_i})^T$. Take a $q\lambda_i \times q\lambda_i$ non-negative matrix $A_i = (a_{lj}^{(i)})$ such that $A_ix_i = \lambda_ix_i$. Then for every $l= 1,...,q\lambda_i$
\begin{equation}\label{equicond}
\sum_{j=1}^{q\lambda_i} a_{lj}^{(i)} = \lambda_i.
\end{equation}
By (\ref{equicond}), exactly $\lambda_i$ edges  start from each vertex in the diagram $B_i$. Clearly, there are several matrices $A_i$ that can  satisfy the above conditions. For instance, we can choose
$$
A_i =
\begin{pmatrix}
2 & 0 & 0 & \ldots & 0 & 0 & 0 &\ldots &\lambda_i - 2\\
0 & 2 & 0 & \ldots & 0 & 0 & 0 &\ldots &\lambda_i - 2\\
\vdots & \vdots & \vdots & \ddots & \vdots & \vdots& \vdots & \ddots & \vdots \\
0 & 0 & 0 & \ldots & 2 & 0 & 0 &\ldots & \lambda_i - 2\\
0 & 0 & 0 & \ldots & 0 & \lambda_i - 1 & 1 & \ldots & 0\\
\vdots & \vdots & \vdots & \vdots &\ddots & \vdots& \vdots & \ddots & \vdots \\
0 & 0 & 0 & \ldots & 0 & 1 & 0 & \ldots & \lambda_i - 1\\
\end{pmatrix}
$$
where the submatrix in the left upper corner with the number two on the main diagonal has the size $i\times i$. Let now $B_i$ be the stationary Bratteli diagram defined by the incidence matrix transpose to $A_i$. The form of $A_i$ means that every minimal component is a 2-odometer and there are exactly $i$ such components. Therefore, the tail equivalence relations $\mathcal R_i$ and $\mathcal R_j$ are not orbit equivalent for diagrams $B_i$ and $B_j$ if $i\neq j$.  Note also that the non-minimal component of $B_i$ has the same form as the diagram $B_0$.

To finish the proof of first part, we conclude that if $\mu_i$ is the measure on $B_i$ defined by the eigenvector $x_i$ and the eigenvalue $\lambda_i$, then $S(\mu_i) = S$ and $\mu_i$ is good by  Corollary \ref{gcd1}. Therefore, $\mu_i$ and $\mu$ are homeomorphic for any $i$.

\textbf{2.} Let $B$ and $\mu$ be as in the theorem and suppose $\mu$ is defined by the eigenvalue $\lambda \in \mathbb{R} \setminus \mathbb{Q}$ of  $A = (a_{ij})_{i,j = 1}^n$, the matrix transposed to the incidence matrix of $B$. To prove the theorem, it suffices to construct a stationary Bratteli diagram $B'$ such that: (i)  there is an ergodic invariant probability good measure $\nu$ on $B'$ such that $S(\nu) = S(\mu)$; (ii)  $B'$ has one more minimal component in comparison with $B$ (in fact, we add another vertex to each level of the initial diagram $B$ and this vertex will determine a minimal component for the tail equivalence relation $\mathcal R'$).

Denote by $(x_1,...,x_n)^T$ the eigenvector corresponding to $\lambda$. Recall that we consider measures on their supports, hence we assume that all $x_i$ are positive.
Let $H = H(x_1,...,x_n)$ denote the additive group generated by $x_1,...,x_n$. Suppose that the vertices $m+1, \ldots, n$ belong to the distinguished class $\alpha$ of vertices  that determine the measure $\mu$. Since $\mu$ is good, there exists $R_0 \in \mathbb{N}$ such that for any integer $R \geq R_0$ the values $\lambda^R x_i$, $i = 1,...,n$, belong to the additive group $H(x_{m+1},...,x_n)$ generated by $\{x_j\}_{j=m+1}^n$ (see Theorem \ref{KritGood}).

Fix $R$ such that $R \geq R_0$.
We will construct a new diagram $B'$ such that the matrix $Q = (q_{ij})_{i,j = 1}^{n+1}$ transposed to the incidence matrix of $B'$ has the  eigenvector
$$
z =  \left(\frac{x_1}{\lambda^R},...,\frac{x_m}{\lambda^R}, \frac{\lambda^R - 1}{\lambda^R}, \frac{x_{m+1}}{\lambda^R},...,\frac{x_n}{\lambda^R}\right)^T
$$
corresponding to the eigenvalue $\psi =  \lambda^M$ where $M = R + N$, $N \in \mathbb{N}$ ($M$ will be chosen below).

To define $Q$, take $A^M = (a_{ij}^{(M)})_{i,j = 1}^n$ and insert  in $A^M$ the additional $(m+1)$-st row $(0,...,0, q_{m+1,m+1}, \ldots, q_{m+1,n + 1})$ and $(m+1)$-st column $(0,...,0, q_{m+1,m+1},0,...0)^T$:
$$
Q =
\begin{pmatrix}
a_{11}^{(M)} & \ldots & a_{1,m}^{(M)} & 0 & a_{1,m+1}^{(M)} & \ldots & a_{1,n}^{(M)}\\
\vdots & \ddots & \vdots  & \vdots & \vdots    & \ddots & \vdots\\
a_{m1}^{(M)} & \ldots & a_{m,m}^{(M)} & 0 & a_{m,m+1}^{(M)} & \ldots & a_{m,n}^{(M)}\\
0 & \ldots & 0 & q_{m+1,m+1} & q_{m+1,m+2} & \ldots & q_{m+1,n + 1}\\
0 & \ldots & 0 & 0 & a_{m+1,m+1}^{(M)} & \ldots & a_{m+1,n}^{(M)}\\
\vdots & \ddots & \vdots & \vdots & \vdots & \ddots & \vdots \\
0 & \ldots & 0 & 0 & a_{n,m+1}^{(M)} & \ldots & a_{n,n}^{(M)}\\
\end{pmatrix}
$$
where $\{q_{m+1,j}\}_{j = m+1}^{n+1}$ are undefined non-negative integers yet. It is worth to mention that some of the entries $a_{ij}^{(M)}$ of $A^M$ may be zero. But the submatrix of $A^M$ formed by the rows enumerated from $m+1$ to $n$ and the columns enumerated from 1 to $m$ is zero matrix since we assumed that the measure $\mu$ is determined by the vertices of the class $\alpha$.

It is obvious that if $q_{m+1,m+1} \geq 2$ and at least one of the numbers $\{q_{m+1,j}\}_{j = m+2}^{n+1}$ is non-zero, then the Bratteli diagram corresponding to $Q$ has one more  minimal component than the diagram $B$ corresponding to $A$. Hence, the dynamical systems $(X_B, \mathcal R)$ cannot be orbit equivalent to the system $(X_{B'}, \mathcal R')$.

Our goal now is to find non-negative integers $\{q_{m+1,j}\}_{j = m+1}^{n+1}$ such that $Q z = \psi z$. It is clear that this equality holds for the the rows $1,..,m,m+2,...n$ of the matrix $Q$ since $A^M\left(\dfrac {x}{\lambda^R}\right) = \lambda^M \dfrac {x}{\lambda^R}$. Therefore, we need only to verify that the equation
\begin{equation}\label{eqirr}
q_{m+1,m+1} \frac{\lambda^R - 1}{\lambda^R} + \sum_{j=m+2}^{n+1} q_{m+1,j} \frac{x_{j-1}}{\lambda^R} = \lambda^N (\lambda^R - 1)
\end{equation}
can be solved for non-negative integers $\{q_{m+1,j}\}_{j = m+1}^{n+1}$.

We will use the geometric representation of algebraic numbers as vectors over $\mathbb Q$ as we did in the proof of Theorem \ref{grouplike}. Let $\lambda$ be the algebraic integer of degree $k$. Suppose $\{\textbf{e}_1,...,\textbf{e}_k\}$ denote the standard basis in $\mathbb R^k$ corresponding to the numbers $1,\lambda,...,\lambda^{k-1}$. Then  (\ref{eqirr}) can be written as follows (we use here and below the notation from the proof of Theorem \ref{grouplike} where, in particular,  matrices $C$ and $D$ were defined):
\begin{equation}\label{eqvec}
\sum_{j=m+2}^{n+1} q_{m+1,j}D^R \textbf{x}_{j-1} = C^N(C^R \textbf{e}_1 - \textbf{e}_1) - q_{m+1,m+1} (\textbf{e}_1 - D^R \textbf{e}_1).
\end{equation}

 Then relation (\ref{eqvec}) can be considered as a $k \times (n - m)$ system of linear equations with respect to the unknowns $\{q_{m+1,j}\}_{j = m+2}^{n+1}$ and the parameter $q_{m+1,m+1}$.
 The matrix $P$ of the system is formed by the columns $\textbf{p}_1,...,\textbf{p}_{n-m}$ where $\textbf{p}_1 = D^R\textbf{x}_{m+1}, ...,\textbf{p}_{n-m} = D^R \textbf{x}_{n}$.
 Since  $\mu$ is good, we see that $\lambda H \subset H \subset H(\frac{x_{m+1}}{\lambda^R},..., \frac{x_{n}}{\lambda^R})$ for $R \geq R_0$. It follows that $\lambda^N (\lambda^R - 1) \in H(\frac{x_{m+1}}{\lambda^R},..., \frac{x_{n}}{\lambda^R})$ for $R \geq R_0$ because $\lambda^R - 1\in H$. Using the correspondence between the elements of $S(\mu)$ and $\mathbb{Q}^k$, we obtain that $C^N(C^R \textbf{e}_1 - \textbf{e}_1) \in H(D^R\textbf{x}_{m+1},...,D^R\textbf{x}_{n})$.
 We show that there exists a natural number $q_{m+1,m+1}$ such that $q_{m+1,m+1} (\textbf{e}_1 - D^R \textbf{e}_1)$ belongs to  $H(D^R\textbf{x}_{m+1},...,D^R\textbf{x}_{n})$. Indeed, by Lemma \ref{finalbasis}, any vector in $\mathbb{Q}^k$ can be represented as a rational linear combination of vectors $D^R\textbf{x}_{m+1},...,D^R\textbf{x}_{n}$. Hence, there exist $t \in \mathbb{N}$ and $\{t_i\}_{i = 1}^{n-m} \subset \mathbb{Z}$ such that $$
 \textbf{e}_1 - D^R \textbf{e}_1 = \sum_{i=1}^{n-m} \frac{t_i}{t} D^R \textbf{x}_{m+i}.
 $$
Therefore, $t (\textbf{e}_1 - D^R \textbf{e}_1) \in H(D^R\textbf{x}_{m+1},...,D^R\textbf{x}_{n})$. It follows that  there exist integers $\{q_{m+1,j}\}_{j = m+2}^{n+1}$ and a non-negative number $q_{m+1,m+1}$ such that relation (\ref{eqvec}) holds.

Let $K = \{\beta_1 D^R \textbf{x}_{m+1} + ... + \beta_{n-m} D^R \textbf{x}_{n} : \beta_1,...,\beta_{n-m} \geq 0 \}$. We consider two cases.

First, let $k = n-m$. Then the  columns of $P$, the vectors $\{D^R \textbf{x}_{m+1},...,D^R \textbf{x}_n\}$,  form a basis in $\mathbb{Q}^k$. Choose $q_{m+1,m+1}$ such that $q_{m+1,m+1} (\textbf{e}_1 - D^R \textbf{e}_1)\in H(D^R\textbf{x}_{m+1},...,D^R\textbf{x}_{n})$. Then the vector in the right part of relation (\ref{eqvec}) has integer coordinates $\{q_{m+1,j}\}_{j = m+2}^{n+1}$ in the basis $\{D^R \textbf{x}_{m+1},...,D^R \textbf{x}_n\}$. We refer now to the proof of Theorem \ref{grouplike} where the behavior of the matrices $C$ and $D$ has been studied. Denote by $l(\textbf{y})$  the line in $\mathbb{R}^k$ generated by a vector $\textbf{y} \in \mathbb{R}^k$. Let $\textbf{y}_1$ be the eigenvector of the matrix $C$ (see the proof of Theorem \ref{grouplike}). We can choose $R$ sufficiently large such that the line $l(\textbf{y}_1)$  lies in $K\bigcup (- K)$. Now we fix $R$ and show that for sufficiently large $N$ the numbers $\{q_{m+1,j}\}_{j = m+2}^{n+1}$ are non-negative. As $N$ tends to infinity, the vector $C^N (C^R \textbf{e}_1 - \textbf{e}_1)$ approaches  to the line $l(\textbf{y}_1)$ in $K$. The norm of this vector tends to infinity.
Hence, for sufficiently large $N$, the right part of  (\ref{eqvec}) lies in $K$ and is an integer combination of linearly independent vectors $\{D^R \textbf{x}_{m+1},...,D^R \textbf{x}_n\}$. Thus, the coefficients $\{q_{m+1,j}\}_{j = m+2}^{n+1}$ of this combination are non-negative integers.

Let $k < n-m$. To find non-negative integer solutions $\{q_{m+1,j}\}_{j = m+2}^{n+1}$ of (\ref{eqvec}), we use a vector generalization of the Frobenius Problem.  The following lemma follows from \cite{Alf, Aliev-Henk}:

\begin{lemm}\label{FrobeniusProblem}
Let $A \in \mathbb{Z}^{m\times n}$, $1 \leq m < n$, be an integral $m\times n$ matrix satisfying

(i) $\gcd(\det(\Omega_i): \Omega_i \mbox{ is an } m\times m \mbox{ minor of }A) = 1$,

(ii) $\{\textbf{y} \in \mathbb{R}^n_{\geq 0} : A\textbf{y} = 0\} = \{0\}$.

\noindent
Denote by $\textbf{v}_1,...,\textbf{v}_n \in \mathbb{Z}^m$  the columns of the matrix $A$, and let
$$
K = \{\beta_1 \textbf{v}_1 + ... + \beta_n \textbf{v}_n : \beta_1,...,\beta_n \geq 0 \}
$$
be the cone generated by $\textbf{v}_1,...,\textbf{v}_n$. Set  $\textbf{v} = \textbf{v}_1 + \cdots + \textbf{v}_n$. Then there exists $0 \leq t_0 < \infty$ such that for $t\geq t_0$ and any vector  $\textbf{b}\in \{t\textbf{v} + K\} \bigcap \mathbb{Z}^m$ there exist a vector  $\textbf{y}$ with non-negative integer entries  such that $A\textbf{y} = \textbf{b}$. (The least possible number  $t_0$ is  called the diagonal Frobenius number $g(A)$).
\end{lemm}

Since the coordinates of  vectors $D^R \textbf{x}_{m+1},...,D^R \textbf{x}_{n}$ in the standard basis $\textbf{e}_1,...,\textbf{e}_k$ are rationals, it is clear that we can multiply the vectors $\textbf{e}_1,...,\textbf{e}_k$ by some rational numbers to obtain a new basis $\{\textbf{e}'_i\}_{i=1}^k$ in which the vectors $\textbf{p}_1,...,\textbf{p}_{n-m}$ have integer coordinates. We will find a basis for which the matrix $P$ written in this basis satisfies conditions (i), (ii) of Lemma \ref{FrobeniusProblem}.

Take the basis $\{\textbf{e}'_i\}$ and denote by $\Lambda = \Lambda (D^R \textbf{x}_{m+1},..., D^R\textbf{x}_n)$ the sublattice of $\mathbb Z^k$ generated by the vectors $\{D^R \textbf{x}_{m+1},..., D^R\textbf{x}_n\}$. In order to have integer solutions for (\ref{eqvec}), the vector $C^N(C^R \textbf{e}_1 - \textbf{e}_1) - q_{m+1,m+1} (\textbf{e}_1 - D^R \textbf{e}_1)$  must belong to $\Lambda$.
The lattice $\Lambda$ has a basis of $k$ elements which belong to $\Lambda$ (see \cite{Nath}). Then we can choose a new basis $\{\textbf{f}_i\}_{i=1}^k$ of $\Lambda$ such that $\Lambda (\textbf{f}_{1},...,\textbf{f}_k)$ is isomorphic to $\mathbb{Z}^k$.  Denote by $P_f$ the matrix $P$ written in the basis $\{\textbf{f}_i\}_{i=1}^k$. Since $D^R \textbf{x}_{m+1},...,D^R\textbf{x}_n \in \Lambda(\textbf{f}_1,...,\textbf{f}_k)$ the matrix $P_f$ has integer entries. We prove that $P_f$ satisfies (i). Consider the identity matrix $I$ (in the basis $\{\textbf{f}_i\}_{i=1}^k$). Since $\Lambda (D^R \textbf{x}_{m+1},...,D^R\textbf{x}_n) = \Lambda(\textbf{f}_1,...,\textbf{f}_k)$, we can express the vectors $\{\textbf{f}_i\}_{i=1}^k$ as integer linear combinations of the vectors $\{D^R \textbf{x}_i\}_{i=1}^k$. Since the determinant is a linear function of its columns, we obtain that $\det I = 1$ is an integer linear combination of the determinants $\det(\Omega_i)$ where $\Omega_i \mbox{ is a } k\times k \mbox{ minor of }P_f$. Hence, condition $(i)$ holds for $P_f$.

It is not hard to prove that condition (ii) holds for $P_f$.
The vectors $\textbf{f}_1,...,\textbf{f}_k$ form the basis of $\mathbb{Q}^k$. Let $J$ be the transition matrix from the basis $\{\textbf{e}_i\}_{i=1}^k$ to the basis $\{\textbf{f}_i\}_{i=1}^k$. Denote by $P_e$ the matrix $P$ in the basis $\{\textbf{e}_i\}_{i=1}^k$. For a vector $\textbf{y} \in \mathbb{Q}^k$ denote by $\textbf{y}_f$ the vector $\textbf y$ written in the basis $\{\textbf{f}_i\}_{i=1}^k$. We keep the notation $\textbf{y}$ for the vector $\textbf{y}$ written in standard basis $\{\textbf{e}_i\}_{i=1}^k$.
The columns of $P_e$ are vectors $D^R \textbf{x}_{m+1},..., D^R\textbf{x}_n$. Since $\textbf{y}_f = J^{-1}\textbf{y}$, we have $P_f = J^{-1}P_e$. Since $\det J \neq 0$, it suffices to prove (ii) for $P_e$. Using the fact that the vectors $D^R \textbf{x}_{m+1},...,D^R\textbf{x}_n$ belong to the same half-space (as shown in the proof of Theorem \ref{grouplike}),  we see that any non-trivial non-negative linear combination of vectors  $D^R \textbf{x}_{m+1},..., D^R\textbf{x}_n$ is non-zero. Indeed, suppose there exist non-negative real numbers $\gamma_1,...,\gamma_{n-m}$ such that $\sum_{i=1}^{n-m} \gamma_i D^R \textbf{x}_{m+i} = 0$. We note that $\langle \sum_{i=1}^{n-m} \gamma_i D^R \textbf{x}_{m+i}, \textbf{n}\rangle = 0$ where $\textbf{n}$ denotes the vector $(1, \lambda, ..., \lambda^{k-1})^T$. On the other hand,
$\langle \sum_{i=1}^{n-m} \gamma_i D^R \textbf{x}_{m+i}, \textbf{n}\rangle = \sum_{i=1}^{n-m} \gamma_i \frac{x_{m+i}}{\lambda^R}$. Since all the values $\frac{x_{m+i}}{\lambda^R}$ are positive, the linear combination is equal to zero if and only if $\gamma_i = 0$ for all $i$.

By Lemma \ref{FrobeniusProblem}, for any right part of (\ref{eqvec}) which belongs to $\{t\textbf{v} + K\} \bigcap \mathbb{Z}^k$ in the basis $\{\textbf{f}_i\}_{i=1}^k$, there exist a non-negative solution $\{q_{m+1,j}\}_{j = m+2}^{n+1}$.
Arguing as in the case $k = n-m$, we show that the vector $C^N(C^R \textbf{e}_1 - \textbf{e}_1) - q_{m+1,m+1} (\textbf{e}_1 - D^R \textbf{e}_1)$ belongs to $\{g(P)\textbf{v} + K\} \bigcap \mathbb{Z}^k$ for sufficiently large $N$. The transformation to the basis $\{\textbf{f}_i\}_{i=1}^k$ alters the coordinates of vectors and the entries of matrices but it doesn't change their properties used in the proof of the case $k = n-m$.

It is left to prove that the measure $\nu$ is good and $S(\nu) = S(\mu)$.
Since  $\mu$ is good, we have that $\lambda^R x_1,...,\lambda^R x_n \in H(x_{m+1},...,x_{n})$ for $R \geq R_0$. It follows from the relation  $1,\lambda^R \in H(x_1,...,x_n)$ that $\lambda^R (\lambda^R - 1) \in H(x_{m+1},...,x_n)$. Therefore,
$$
\lambda^R \dfrac{x_1}{\lambda^R},...,\lambda^R \dfrac{x_n}{\lambda^R},\lambda^R \dfrac{\lambda^R - 1}{\lambda^R} \in H\left(\dfrac{x_{m+1}}{\lambda^R},...,\dfrac{x_{n}}{\lambda^R}\right)
$$
for $R \geq R_0$.
Hence, we obtain that $\lambda^{R} z_1,...,\lambda^{R} z_{n+1} \in H(z_{m+1},...,z_{n+1})$. Since $M > R$ we have $\psi z_1,...,\psi z_{n+1} \in H(z_{m+1},...,z_{n+1})$ and $\nu$ is good by Theorem \ref{KritGood}.

Finally, we conclude that
 $$
 H(z_1,...,z_{n+1}) = \dfrac{1}{\lambda^{R}} H(x_{1},...,x_{n}),\ \ \  \dfrac{1}{\psi} H(z_1,...,z_{n+1}) = \dfrac{1}{\lambda^{R+M}} H(x_{1},...,x_{n}).
 $$
Since $S(\mu) = \bigcup_{N \in \mathbb{N}} \frac{1}{\lambda^N}H(x_1,...,x_n)$ and $\lambda H \subset H$, we have $S(\nu) = S(\mu)$.
\hfill$\blacksquare$
\medskip

\begin{remar}
For $\lambda \notin \mathbb Q$, we can proceed as in the rational case by finding a measure $\nu$ on a simple stationary Bratteli diagram such that  $S(\nu) = S(\mu)$. To do this, we construct a diagram with an $n\times n$ matrix $\widetilde{A} = (\widetilde{a}_{ij})$ transpose to the incidence matrix such that $\widetilde{A}x = \lambda^M x$ for some $M \in \mathbb{N}$. We obtain the matrix $\widetilde{A}$ by taking $A^M$ for sufficiently large $M$ and changing the zero block $\{a_{ij}^{(M)}\}_{i = m+1, j = 1}^{n,m}$ to a non-zero one.
We have $q_i \textbf{x}_i = \sum_{j = m+1}^n p_{ij} \textbf{x}_j$ for $i = 1,...,m$ and some $q_i \in \mathbb{N},\ p_{ij} \in \mathbb{Z}$. Since $A^M x = \lambda^M x$, we have $\sum_{i=m+1}^n a^{(M)}_{ij}x_j = \lambda^M x_i$ for $i = m+1,...,n$. The block $\{a_{ij}\}_{i,j = m+1}^{n}$ of the matrix $A$ is positive, hence we can make the numbers $\{a_{ij}^{(M)}\}_{i,j = m+1}^{n}$ arbitrary large. In particular, we take $M$ such that $a_{m+1,j}^{(M)} \geq p_{1,j}$ for all $j = m+1,...,n$. Then $q_1 x_1 + \sum_{j=m+1}^n (a_{m+1,j}^{(M)} - p_{1,j})x_j = \sum_{j=m+1}^n a_{m+1,j}^{(M)}x_j = \lambda x_{m+1}$ and $\widetilde{a}_{m+1,1} = q_1 > 0$. We proceed with other elements of the zero block in the similar way.

Thus, we obtain that for any finite ergodic invariant measure $\mu$ on a stationary Bratteli diagram there exists a good measure $\nu \in \mathcal S$ such that $S(\mu) = S(\nu)$. This result together with an example is discussed in the last section.
\end{remar}

We will also construct two measures $\mu_1$ and $\mu_2$ such that $S(\mu_1) = S(\mu_2)$, the measure $\mu_1$ is good and the measure $\mu_2$ is not good. Hence these measures are not homeomorphic. The example can be found in the next section.

Now we consider two classes of measures on stationary Bratteli diagrams: (i) measures of Bernoulli type, (ii) measures satisfying the Quotient Condition.

\begin{prop}\label{Qcond}
Let $\mu$ be an ergodic $\mathcal R$-invariant measure on a stationary Bratteli diagram $B$. Then $\mu$ does not satisfy the Quotient Condition.
\end{prop}

\noindent
\textbf{Proof.}
Let $A = F^T$ be the $n \times n$ matrix transposed to the incidence matrix of the diagram $B$. Let $\lambda$ and $x = (x_1,...,x_n)^T$ be the eigenvalue and eigenvector of $A$ which generate $\mu$. Denote by  $H$ the additive subgroup of $\mathbb{R}$ generated by $\{x_1, \ldots , x_n\}$.
Let eigenvalue $\lambda$ be the algebraic number of degree $k$. If $k = 1$, the clopen values set $S(\mu)$ consists of rational numbers. To get a contradiction, we assume that $\mu$ satisfies the Quotient Condition. Then, by Theorems \ref{BernQuot} and \ref{grouplike}, we have
$$
\mathbb{Q} \bigcap [0,1]\subset S(\mu) = \left(\bigcup_{N=0}^\infty \frac{1}{\lambda^N} H\right) \bigcap [0,1].
$$
Then for any $l \in \mathbb{N}$ there exists $N \in \mathbb{N}$ such that $\frac{1}{l} \in \frac{1}{\lambda^N}H$. In the vector form, the latter is written as $\frac{1}{l} \textbf{e}_1 \in H(D^N \textbf{x}_1,...,D^N \textbf{x}_n)$. The matrix $D$ (used in the proof of Theorem \ref{grouplike}) has rational entries $(d_{ij})$. There exist non-negative integer $t$ such that  $d_{ij} = \frac{t_{ij}}{t}$ where $t_{ij} \in \mathbb{Z}$ for $i = 1,...,k,\ j=1,...,n$. Denote by $T$ the $k\times k$ matrix with entries $(t_{ij})$.  Then $\frac{1}{l} \textbf{e}_1 \in \frac{1}{t^N} H(T^N \textbf{x}_1,...,T^N \textbf{x}_n)$. The set $L$ of prime divisors of $t$ and the denominators of the coordinates of $\textbf{x}_1,...,\textbf{x}_n$ is finite. We can choose $l$ as a prime number which does not belong to $L$. Then $\frac{1}{l} \textbf{e}_1$ cannot belong to $H(D^N \textbf{x}_1,...,D^N \textbf{x}_n)$ for any $N \in \mathbb{N}$. This is a contradiction. \hfill$\blacksquare$
\medskip

Now we focus on  rational measures, i.e. on the case when the clopen values sets belong to $\mathbb Q$.

\begin{prop}\label{goodBern}
Let $\mu$ be a rational measure on a stationary diagram $B$ defined by a  distinguished eigenvalue $\lambda$ of the matrix $A = F^T$. Denote by $(\frac{p_1}{q}, \ldots, \frac{p_n}{q})$  the corresponding reduced vector.

(1) A number $y \in S(\mu) + \mathbb{Z}$ if and only if there exists $N \in \mathbb N$ such that $y \lambda^{N} q  \in \mathbb{Z}$.

(2) The set $S(\mu)$ is multiplicative if and only if there exists $K \in \mathbb{N}$ such that $q \mid \lambda^{K}$.

(3) If $\mu$ is good, then $\mu$ is of Bernoulli type if and only if there exists $K \in \mathbb{N}$ such that $q \mid \lambda^{K}$.

\end{prop}

\noindent
\textbf{Proof.}
\textbf{1.} Suppose there exists $N \in \mathbb{N}$ such that $y \lambda^{N} q \in \mathbb{Z}$. Since $\gcd(p_{1},\ldots,p_{n}) = 1$, there exist $u_{1}, \ldots , u_{n} \in \mathbb{Z}$ such that $u_{1}p_{1} + \ldots + u_{n}p_{n} = 1$. Then one can take $k_{1}, \ldots , k_{n} \in \mathbb{Z}$ such that $y \lambda^{N} q = k_{1}p_{1} + \ldots + k_{n}p_{n}$. Hence $y = \frac{1}{\lambda^{N}}(k_{1}\frac{p_{1}}{q} + \ldots + k_{n}\frac{p_{n}}{q})$. By Theorem \ref{grouplike}, $S(\mu)$ is group-like. Therefore $y \in S(\mu) + \mathbb{Z}$.

Conversely, let $y \in S(\mu) + \mathbb{Z}$. Then there exist $k_{1}, \ldots , k_{n} \in \mathbb{Z}$ and $N \in \mathbb{N}$ such that $y = \frac{1}{\lambda^{N}}(k_{1}\frac{p_{1}}{q} + \ldots + k_{n}\frac{p_{n}}{q})$. Hence $y \lambda^{N} q = k_{1}p_{1} + \ldots + k_{n}p_{n} \in \mathbb{Z}$.

\textbf{2.} Suppose $S(\mu)$ is multiplicative. Take  $u_{1}, \ldots , u_{n} \in \mathbb{Z}$ such that $u_{1}p_{1} + \ldots + u_{n}p_{n} = 1$. Since $S(\mu)$ is group-like, the fraction $\frac{1}{q} = u_{1}\frac{p_{1}}{q} + \ldots + u_{n}\frac{p_{n}}{q}$ lies in $S(\mu)$. Then the fraction $\frac{1}{q^2}$ belongs to $S(\mu)$. Hence, $\frac{1}{q^{2}} \lambda^{K} q = \frac{1}{q} \lambda^{K} \in \mathbb{N}$ for some $K \in \mathbb{N}$.

Suppose $q \mid \lambda^{K}$. Let $y_{1}, y_{2} \in S(\mu)$. Then $y_{1} \lambda^{M} q$ and $y_{2} \lambda^{N} q$ are integers for some $M, N \in \mathbb{N}$. Hence, $y_{1} y_{2} \lambda^{M + N} q^{2} \in \mathbb{Z}$. Then $y_{1} y_{2} \lambda^{M + N + K} q$ is an integer. Therefore, $y_{1} y_{2} \in S(\mu)$.

\textbf{3.} Let $q \mid \lambda^{K}$ for some $K \in \mathbb{N}$. Then the set $S(\mu)$ is multiplicative, and $S(\mu)$ is ring-like. By Theorem \ref{BernQuot}, it suffices to show that every positive element $y$ from  $S(\mu) + \mathbb{Z}$ is a sum of positive units of $S(\mu) + \mathbb{Z}$. There exist $N \in \mathbb{N}$ and $k_{1}, \ldots, k_{n}$ such that $y =  \frac{1}{\lambda^{N} q}(k_{1}p_{1} + \ldots + k_{n}p_{n})$. The fractions $\frac{1}{\lambda^{N}}$ and $\frac{1}{q}$ are the positive units of the ring $S(\mu) + \mathbb{Z}$. Therefore, their product is a unit of the ring. Since $y >0 $,  the integer $k_{1}p_{1} + \ldots + k_{n}p_{n}$ is non-negative. Thus, $y$ is the sum of the positive unit $\frac{1}{q \lambda^{N}}$ taken  $k_{1}p_{1} + \ldots + k_{n}p_{n}$ times.

Conversely, let $\mu$ be a good measure of Bernoulli type. Then $S(\mu)$ is ringlike. Hence, there exists $K \in \mathbb{N}$ such that $q \mid \lambda^{K}$.  \hfill$\blacksquare$


\section{Examples}

In this section, we consider several examples of ergodic invariant measures on some stationary Bratteli diagrams illustrating the results proved in Section 3.
\medskip

\textbf{Example 1.} Let $B$ be the stationary Bratteli diagram  with incidence matrix
$$
F =
\begin{pmatrix}
1 & 1 & 0\\
1 & 2 & 0\\
0 & 1 & 3 \\
\end{pmatrix}.
$$
This diagram $B$ looks as follows:
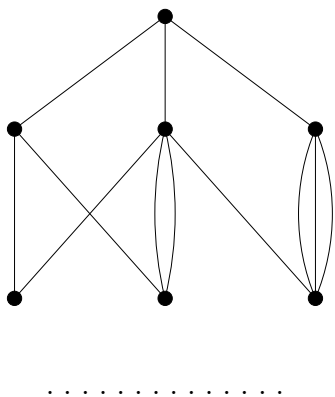
\unitlength = 0.5cm
\begin{center}
\begin{graph}(8,13)
\graphnodesize{0.4}
%
\roundnode{V0}(5,12)
\roundnode{V11}(1,9)
\roundnode{V12}(5,9)
\roundnode{V13}(9,9)
\edge{V11}{V0}
\edge{V12}{V0}
\edge{V13}{V0}
\roundnode{V21}(1,4.5)
\roundnode{V22}(5,4.5)
\roundnode{V23}(9,4.5)
\edge{V21}{V11}
\edge{V21}{V12}
\bow{V22}{V12}{0.06}
\bow{V22}{V12}{-0.06}
\edge{V22}{V11}
\bow{V23}{V13}{0.10}
\bow{V23}{V13}{-0.10}
\edge{V23}{V13}
\edge{V23}{V12}
\freetext(5,2){$.\ .\ .\ .\ .\ .\ .\ .\ .\ .\ .\ .\ .\ .$}
\end{graph}
\end{center}

This diagram has two simple stationary subdiagrams sitting on the first two and the third vertices, respectively. The left subdiagram is a minimal component for $\mathcal R$.

We consider the  $\mathcal R$-invariant ergodic measures for the diagram $B$ and compute their clopen values sets and prove some properties of these  measures.

Denote $A = F^T$. The eigenvectors $x = (\frac{3 - \sqrt{5}}{2}, \frac{\sqrt{5} - 1}{2}, 0)^T$ and $y = (\frac{1}{4}, \frac{1}{2}, \frac{1}{4})^T$ of the matrix $A$ correspond to the eigenvalues $\lambda_{1} = \frac{3 + \sqrt{5}}{2}$ and $\lambda_{2} = 3$, respectively.
Thus, there are two probability ergodic $\mathcal{R}$-invariant measures on the path space of the diagram $B$. Denote by $\mu_{1}$ the measure generated by the vector $x$ and the eigenvalue $\lambda_{1}$, and by $\mu_{2}$ the measure generated by the vector $y$ and the eigenvalue $\lambda_{2}$. We note that $S(\mu_1) \cap (\mathbb R \setminus \mathbb Q) \neq \emptyset$ and $S(\mu_2) \subset \mathbb Q$.

Let  $h^{(n)} = (h^{(n)}_{i})$ denote the heights of  towers  corresponding to the vertices $i = 1, 2, 3$ enumerated from left to right. We have  $h^{(1)} = (1, 1, 1)$ and
\begin{equation}\label{hcoordin}
h^{(n + 1)} = Fh^{(n)}.
\end{equation}
In order to find the solutions $h^{(n)}$ of (\ref{hcoordin}), we use generating functions. Set $f^{(i)}(s) = \sum_{n = 0}^{\infty} h_{i}^{(n + 1)}s^{n}$.
It can be shown  that
\begin{eqnarray*}
f^{(1)}(s) & = & \frac{1 - s}{1 - 3s + s^{2}},\\
f^{(2)}(s) & = & \frac{1}{1 - 3s + s^{2}},\\
f^{(3)}(s) & = &\frac{(1 - s)^{2}}{(1 - 3s + s^{2})(1 - 3s)}.
\end{eqnarray*}
Decomposing the generating functions into the series,  we obtain
\begin{eqnarray*}
h_{1}^{(n)} & = & \frac{1}{\sqrt{5}}\left(\frac{1 + \sqrt{5}}{2} \left(\frac{3 + \sqrt{5}}{2}\right)^{n - 1} - \frac{1 - \sqrt{5}}{2}\left(\frac{3 - \sqrt{5}}{2}\right)^{n-1} \right),\\
h_{2}^{(n)} & = & \frac{1}{\sqrt{5}}\left( \left(\frac{3 + \sqrt{5}}{2}\right)^{n} - \left(\frac{3 - \sqrt{5}}{2}\right)^{n} \right),\\
h_{3}^{(n)} & = & 4 \cdot 3^{n} + \frac{7\sqrt{5} - 15}{10} \left(\frac{3 - \sqrt{5}}{2}\right)^{n-1} - \frac{7\sqrt{5} + 15}{10}\left(\frac{3 + \sqrt{5}}{2}\right)^{n-1}.
\end{eqnarray*}
We see that $h^{(n)}_1 = f_{2n - 1}$ and $h^{(n)}_2 = f_{2n}$ where $f_i$ is the $i$-th Fibonacci number.

The above computation allows one to determine explicitly all elements of the sets $S(\mu_{1})$ and $S(\mu_{2})$.
The ergodic measures $\mu_{1}$ and $\mu_{2}$ are good. Indeed, the measure $\mu_1$ is supported on a simple subdiagram, hence $\mu_1$ is good. It follows from Theorem \ref{KritGood} that $\mu_2$ is also good. One can show (see Example 2) that in the simplex of all $\mathcal R$-invariant probability measures only these ergodic measures are good. In other words, the measure $\nu_\alpha = \alpha\mu_1 + (1- \alpha)\mu_2$ is not good for any $\alpha \in (0,1)$. We remark also that, by Proposition \ref{Qcond}, neither $\mu_1$ nor $\mu_2$ satisfy the Quotient Condition.

\begin{prop} Let $B, \mu_1, \mu_2$ be as above. Then $\mu_{1}$ is of Bernoulli type and $\mu_{2}$ is not.
\end{prop}

\noindent
\textbf{Proof}. We use Theorem \ref{BernQuot} to prove that $\mu_1$ is of Bernoulli type.
Since $\frac{2}{3 + \sqrt{5}} = \frac{3 - \sqrt{5}}{2}$, we see that
\begin{equation}\label{Smu1}
S(\mu_{1}) = \left\{k^{(n)}_{1}\left(\frac{3 - \sqrt{5}}{2}\right)^{n} + \; k^{(n)}_{2}\frac{\sqrt{5} - 1}{2}\left(\frac{3 - \sqrt{5}}{2}\right)^{n-1}: n \in \mathbb{N}\right\}
\end{equation}
where $0\leq k^{(n)}_{i} \leq h^{(n)}_{i}$.

Denote $G = S(\mu) + \mathbb{Z}$. By Theorem \ref{goodmeasure}, $G$ is an additive subgroup of $\mathbb{R}$. Since $\left(\frac{\sqrt{5} - 1}{2}\right)^{2} = \frac{3 - \sqrt{5}}{2}$, the group $G$ is multiplicative. Hence $G$ is a ring. We note that the number $1$ is a positive unit of the ring $G$, and the fractions $\frac{\sqrt{5} - 1}{2}$ and $\frac{3 - \sqrt{5}}{2}$ are also positive units of  $G$ because $\frac{2}{\sqrt{5} - 1} = \frac{\sqrt{5} - 1}{2} + 1 \in G$ and $\frac{2}{3 - \sqrt{5}} = 3 - \frac{3 - \sqrt{5}}{2} \in G$.  Therefore, every positive element of $G$ is a sum of positive units of $G$. By Theorem \ref{BernQuot}, the measure $\mu_{1}$ is of Bernoulli type.
\medskip

We have already observed that the measure $\mu_{2}$ is good. Then $\mu_{2}$ is not of Bernoulli type by Proposition \ref{goodBern}.
\hfill$\blacksquare$
\\

\textbf{Example 2.} The following example is a generalization of Example 1. We consider a class of non-simple stationary Bratteli diagrams that have one minimal and one non-minimal component and have exactly two ergodic probability $\mathcal R$-invariant measures defined  by these components. Let $B$ be such a diagram. Suppose that $S(\mu_1)$ contains irrational numbers for the measure $\mu_1$ supported on the minimal component, and $S(\mu_2) \subset \mathbb Q$ for the other ergodic measure  $\mu_2$. The Bratteli diagram in Example 1 belongs to this class. Since $\mu_1$ and $\mu_2$ are the only ergodic measures, any $\mathcal R$-invariant measure $\nu$ is of the form $\nu_{\alpha} = \alpha \mu_{1} + (1 - \alpha)\mu_{2}$ where $\alpha \in [0,1]$. These measures form the convex simplex of all $\mathcal{R}$-invariant probability measures on $X_{B}$. Our goal is to show that the measure $\nu_{\alpha}$ is good only for $\alpha = 0,1$.

\begin{prop}\label{simplex} Let $\mu_1, \mu_2$ and $\nu_\alpha$ be as above. Then the measure
$\nu_{\alpha}$ on the Cantor space $X_{B}$ is not good for any $\alpha \in (0,1)$;
\end{prop}

\noindent
\textbf{Proof.} Let $A = F^T$ be the $n \times n$ matrix transposed to the incidence matrix of the diagram $B$. Assume the vertices $1,...,m$ belong to the minimal component of the diagram and the vertices $m+1,...,n$ are in the non-minimal component.
Let $\mu_1$ be generated by the eigenvalue $\lambda_1 \in \mathbb{R}\setminus \mathbb{Q}$ and the eigenvector $x = (x_1,...,x_m,0,...,0)^T$, and  let $\mu_2$ be generated by the eigenvalue $\lambda_2 \in \mathbb{N}$ and the eigenvector $y = (\frac{p_1}{q},...,\frac{p_n}{q})^T$ where $p_i, q \in \mathbb{N}$.

We first chose two particular clopen sets $U$ and $V$ such that $\nu_{\alpha}(U) < \nu_{\alpha}(V)$ and then we show that there is no clopen subset of $V$ with measure $\nu_{\alpha}(U)$. Denote by $V$ the cylinder set of length 1 that ends in a vertex of the non-minimal component. Without loss of generality, we assume that $V$ ends in the vertex number $n$. Denote by $U_N$ the cylinder set of length $N$ which ends in the vertex $1$ of the minimal component of $B$.
Take $N$  sufficiently large so that $\nu_{\alpha}(U_{N}) < \nu_{\alpha}(V)$. Let $U = U_{N}$.

Suppose $\nu_{\alpha}$ is good for some $\alpha \in (0,1)$. We have
$$
\nu_{\alpha}(U) = \alpha \frac{x_1}{\lambda_1^{N-1}} + (1 - \alpha)\frac{p_1}{q\lambda_2^{N-1}} < \nu_{\alpha}(V) = \frac{p_n}{q}.
$$
By assumption, there exists a subset $W \subset V$ such that $\nu_{\alpha}(U) = \nu_{\alpha}(W)$. Since $W \subset V$, there exist integers $k$, $M$ (we can always choose $M > N$) such that
$$
\nu_{\alpha}(W) = (1 - \alpha) \frac{k}{q \lambda_2^{M}}.
$$
It follows from the equality $\nu_{\alpha}(W) = \nu_{\alpha}(U)$ that
\begin{equation}\label{k}
k = \frac{\alpha}{1 - \alpha} \cdot \frac{x_1 q \lambda_2^M}{\lambda_1^{N-1}} + p_1\lambda_2^{M - N + 1}.
\end{equation}
Since the numbers $k,\ p_1\lambda_2^{M - N + 1},\ q\lambda_2^M$ are integers, we have
$$
\frac{\alpha}{1 - \alpha} \cdot \frac{x_1}{\lambda_1^{N-1}} \in \mathbb{Q}.
$$
We can repeat the same arguments for $N+1$ instead of $N$ and obtain that
$$
\frac{\alpha}{1 - \alpha} \cdot \frac{x_1}{\lambda_1^{N}} \in \mathbb{Q}.
$$
Hence, the ratio of the two above mentioned values should be rational. But this ratio equals $\lambda_1 \in \mathbb{R} \setminus \mathbb{Q}$. This is a contradiction.
\\

\textbf{Example 3.} We consider now a class of stationary Bratteli diagrams and determine which measures on them are good. Fix an integer  $N \geq 3$ and let
$$
F_N =
\begin{pmatrix}
2 & 0 & 0\\
1 & N & 1\\
1 & 1 & N \\
\end{pmatrix}
$$
be the incidence matrix of the Bratteli diagram $B_N$. For $A_N = F_N^T$ we easily find the Perron-Frobenius eigenvalue $\lambda = N+1$ and the corresponding probability eigenvector
$$
x = \left(\frac{1}{N},\ \frac{N-1}{2N},\ \frac{N-1}{2N}\right)^T.
$$
Let $\mu_N$ be the measure on $B_N$ determined by $\lambda$ and the eigenvector  $x$. It follows from Theorem \ref{KritGood} that  $\mu_N$ is a good measure if and only if  for all sufficiently large $R$ we have
$$
\frac{(N+1)^R}{N} \in  \frac{N-1}{2N}\mathbb Z
$$
or, equivalently, $\frac{2(N+1)^R}{N-1}$ is an integer, $k\in \mathbb N$. This is possible if and only if $N= 2^k+1$. For instance, the measure $\mu_N$ is good for $N = 3, 5$ but is not good for $N= 4$.
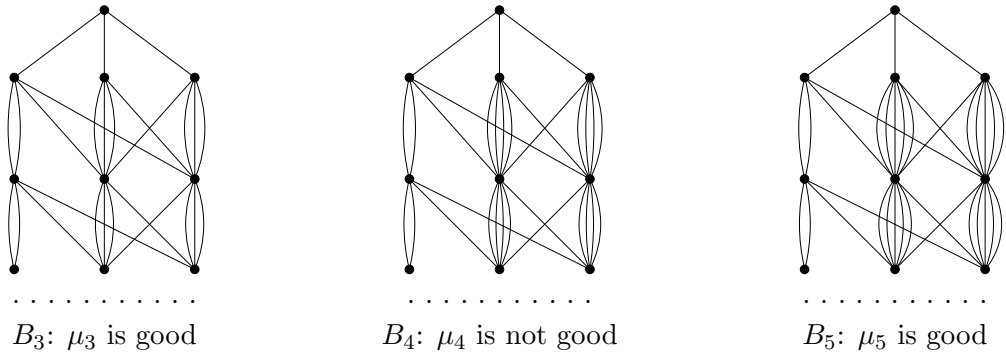
\unitlength = 0.3cm
\begin{center}
\begin{tabular*}{0.99\textwidth}%
     {@{\extracolsep{\fill}}ccc}

\begin{graph}(10,13)
\graphnodesize{0.4}
\roundnode{V0}(5,12)
\roundnode{V11}(1,9)
\roundnode{V12}(5,9)
\roundnode{V13}(9,9)
\edge{V11}{V0}
\edge{V12}{V0}
\edge{V13}{V0}
\roundnode{V21}(1,4.5)
\roundnode{V22}(5,4.5)
\roundnode{V23}(9,4.5)
\bow{V21}{V11}{0.06}
\bow{V21}{V11}{-0.06}
\edge{V22}{V11}
\bow{V22}{V12}{0.10}
\edge{V22}{V12}
\bow{V22}{V12}{-0.10}
\edge{V22}{V13}
\bow{V23}{V13}{0.10}
\edge{V23}{V13}
\bow{V23}{V13}{-0.10}
\edge{V23}{V11}
\edge{V23}{V12}
\roundnode{V31}(1,0.5)
\roundnode{V32}(5,0.5)
\roundnode{V33}(9,0.5)
\bow{V31}{V21}{0.06}
\bow{V31}{V21}{-0.06}
\edge{V32}{V21}
\bow{V32}{V22}{0.10}
\edge{V32}{V22}
\bow{V32}{V22}{-0.10}
\edge{V32}{V23}
\bow{V33}{V23}{0.10}
\edge{V33}{V23}
\bow{V33}{V23}{-0.10}
\edge{V33}{V21}
\edge{V33}{V22}
\freetext(5,-.9){$.\ .\ .\ .\ .\ .\ .\ .\ .\ .\ .$}%
 \freetext(5,-2.5){{$B_3$: $\mu_3$ is good}}
\end{graph}

            &
\unitlength = 0.3cm
\begin{graph}(10,13)
\graphnodesize{0.4}
\roundnode{V0}(5,12)
\roundnode{V11}(1,9)
\roundnode{V12}(5,9)
\roundnode{V13}(9,9)
\edge{V11}{V0}
\edge{V12}{V0}
\edge{V13}{V0}
\roundnode{V21}(1,4.5)
\roundnode{V22}(5,4.5)
\roundnode{V23}(9,4.5)
\bow{V21}{V11}{0.06}
\bow{V21}{V11}{-0.06}
\edge{V22}{V11}
\bow{V22}{V12}{0.04}
\bow{V22}{V12}{0.12}
\bow{V22}{V12}{-0.12}
\bow{V22}{V12}{-0.04}
\edge{V22}{V13}
\bow{V23}{V13}{0.04}
\bow{V23}{V13}{0.12}
\bow{V23}{V13}{-0.12}
\bow{V23}{V13}{-0.04}
\edge{V23}{V11}
\edge{V23}{V12}
\roundnode{V31}(1,0.5)
\roundnode{V32}(5,0.5)
\roundnode{V33}(9,0.5)
\bow{V31}{V21}{0.06}
\bow{V31}{V21}{-0.06}
\edge{V32}{V21}
\bow{V32}{V22}{0.04}
\bow{V32}{V22}{0.12}
\bow{V32}{V22}{-0.12}
\bow{V32}{V22}{-0.04}
\edge{V32}{V23}
\bow{V33}{V23}{0.04}
\bow{V33}{V23}{0.12}
\bow{V33}{V23}{-0.12}
\bow{V33}{V23}{-0.04}
\edge{V33}{V21}
\edge{V33}{V22}
\freetext(5,-.9){$.\ .\ .\ .\ .\ .\ .\ .\ .\ .\ .$}%
 \freetext(5,-2.5){{$B_4$: $\mu_4$ is not good}}
\end{graph}
            &
\unitlength = 0.3cm
\begin{graph}(10,13)
\graphnodesize{0.4}
\roundnode{V0}(5,12)
\roundnode{V11}(1,9)
\roundnode{V12}(5,9)
\roundnode{V13}(9,9)
\edge{V11}{V0}
\edge{V12}{V0}
\edge{V13}{V0}
\roundnode{V21}(1,4.5)
\roundnode{V22}(5,4.5)
\roundnode{V23}(9,4.5)
\bow{V21}{V11}{0.06}
\bow{V21}{V11}{-0.06}
\edge{V22}{V11}
\bow{V22}{V12}{0.09}
\bow{V22}{V12}{0.18}
\edge{V22}{V12}
\bow{V22}{V12}{-0.18}
\bow{V22}{V12}{-0.09}
\edge{V22}{V13}
\bow{V23}{V13}{0.09}
\bow{V23}{V13}{0.18}
\edge{V23}{V13}
\bow{V23}{V13}{-0.18}
\bow{V23}{V13}{-0.09}
\edge{V23}{V11}
\edge{V23}{V12}
\roundnode{V31}(1,0.5)
\roundnode{V32}(5,0.5)
\roundnode{V33}(9,0.5)
\bow{V31}{V21}{0.06}
\bow{V31}{V21}{-0.06}
\edge{V32}{V21}
\bow{V32}{V22}{0.09}
\bow{V32}{V22}{0.18}
\edge{V32}{V22}
\bow{V32}{V22}{-0.18}
\bow{V32}{V22}{-0.09}
\edge{V32}{V23}
\bow{V33}{V23}{0.09}
\bow{V33}{V23}{0.18}
\edge{V33}{V23}
\bow{V33}{V23}{-0.18}
\bow{V33}{V23}{-0.09}
\edge{V33}{V21}
\edge{V33}{V22}
\freetext(5,-.9){$.\ .\ .\ .\ .\ .\ .\ .\ .\ .\ .$}%
 \freetext(5,-2.5){{$B_5$: $\mu_5$ is good}}
\end{graph}
\vspace{1.5cm}
\end{tabular*}
\end{center}

In the case when $N = 4$,
we see that the eigenvector corresponding to $\lambda = 5$  is $(\frac{2}{8}, \frac{3}{8}, \frac{3}{8})$.  The fact that $\mu_4$ is not good can be also proved straightforwardly. Indeed, the measure of any cylinder set that ends in one of the last two vertices of the diagram $B$ at the level $n$ is $\frac{3}{8 \cdot 5^{n-1}}$. Hence, the measure $\mu_4(V)$ of any clopen set $V$ combined from such cylinder sets is a rational number with factor $3$ in the numerator. Let, for instance, $U$ be the cylinder set of length 1 that ends at the first vertex of the diagram and $V$ be the cylinder set of length 1 that ends at the second vertex. Then  $\frac{2}{8} = \mu_4(U) < \mu_4(V) = \frac{3}{8}$ and the measure of any clopen subset of $V$ is a rational irreducible fraction with the number $3$ as a numerator factor. Thus, there is no clopen subset of $V$ with measure $\mu_4(U) = \frac{2}{8}$ and $\mu_4$ is not good.
\medskip

On the other hand, it is easy to find a  measure $\nu$ on a simple stationary Bratteli diagram such that $S(\mu_4) = S(\nu)$. In fact, the following general statement holds.

\begin{prop}
For any probability ergodic $\mathcal R$-invariant measure $\mu$ on a stationary Bratteli diagram there exists a good measure $\nu \in \mathcal S$ such that $S(\mu) = S(\nu)$.
\end{prop}

\noindent
\textbf{Proof}. The proof immediately follows from Theorem \ref{goodhomeo} and the remark after the theorem: there exists a measure $\nu$ on a simple Bratteli diagram (hence $\nu$ is good) such that $S(\mu) = S(\nu)$. \hfill$\blacksquare$
\medskip

Another way to find such a measure $\nu$ is to find a non-negative integer solution $(p_{ij}) = P$ of the system $Py = \psi y$ where the eigenvector $y$ and eigenvalue $\psi$ generate the measure $\nu$ and are chosen such that $S(\nu) = S(\mu)$. For instance, if $\mu = \mu_4$, then we can take
$$
P =
\begin{pmatrix}
1 & 2 & 0\\
1 & 2 & 1\\
9 & 3 & 2 \\
\end{pmatrix}
$$
where $y = (\frac{1}{8}, \frac{2}{8}, \frac{5}{8})$ and $\psi = 5$. It is not hard to see that the clopen values sets for  $S(\nu)$ and $S(\mu_4)$ coincide, and the measure $\nu$ is good but  $\mu_4$ is not. Hence, by Theorem \ref{goodmeasure}, the measures $\mu_4$ and $\nu$ are not topologically equivalent.


\begin{thebibliography}{99}

\bibitem{Akin1}
E.\,Akin, Measures on Cantor space, ``Topology Proc.'',  \textbf{24} (1999),  1 - 34.

\bibitem{Akin2}
E.\,Akin, Good measures on Cantor space, ``Trans. Amer. Math. Soc.'', \textbf{357} (2005),  2681 - 2722.

\bibitem{Akin3}
E. Akin, R. Dougherty, R.D. Mauldin, A. Yingst, Which Bernoulli measures are good measures? ``Colloq. Math.'', \textbf{110} (2008),  243 -- 291.


\bibitem{Alf}
J.\,L.\,Ramirez Alfonsin, The Diophantine Frobenius Problem, Oxford Lecture Series in Mathematics and its Applications, 30, Oxford University Press, 2005.

\bibitem{Aliev-Henk}
I.\,Aliev, M.\,Henk, On feasibility of integer knapsacks, 	 arXiv:0911.4186v1

\bibitem{Alp-Pr}
S.\,Alpern and V.\,S.\,Prasad, Typical Dynamics of Volume Preserving Homeomorphisms, Cambridge Tracts in Mathematics, 139, Cambridge Univercity Press, Cambridge 2000.

\bibitem{Austin}
T.\,D.\,Austin, A pair of non-homeomorphic product measures on the Cantor set, ``Math. Proc. Cam. Phil. Soc.'', \textbf{142} (2007), 103 - 110.

\bibitem{S.B.}
S. Bezuglyi, J.Kwiatkowski, K.Medynets, and B.Solomyak, Invariant measures on stationary Bratteli diagrams, ``Ergodic Theory Dynam. Syst.'', \textbf{30}  (2010), 973 - 1007.

\bibitem{BKM}
S.~Bezuglyi, J.~Kwiatkowski, and K.~Medynets,  Aperiodic substitution  systems and their Bratteli diagrams,  ``Ergodic Theory Dynam. Syst.'', \textbf{29} (2009), 37 - 72.


\bibitem{D-M-Y}
R.\,Dougherty, R.\,Daniel Mauldin, and A.\,Yingst, On homeomorphic Bernoulli measures on the Cantor space, ``Trans. Amer. Math. Soc.'' \textbf{359} (200 ), 6155 - 6166

\bibitem{GW}
E. Glasner, B. Weiss, Weak orbital equivalence of minimal Cantor systems, ``Internat. J. Math.'' \textbf{6} (1995), 559 - 579.

\bibitem{Nath}
Melvyn B. Nathanson, Additive Number Theory: Inverse Problems and the Geometry of Sumsets, Springer-Verlag, New York, 1996.

\bibitem{Navarro}
F.J. Navarro-Bermudez, Topologically equivalent measures in the Cantor space, ``Proc. Amer. Math. Soc.'',  \textbf{77} (1979), 229 - 236.

\bibitem{Navarro-Oxtoby}
F.J. Navarro-Bermudez and J.C. Oxtoby, Four topologically equivalent measures in the Cantor space, ``Proc. Amer. Math. Soc.'', \textbf{104} (1988), 229 - 236.

\bibitem{Oxt-Pr}
J.\,C.\,Oxtoby and V.\,S.\,Prasad, Homeomorphic measures in the Hilbert Cube, ``Pacific J. Math.'', \textbf{77} (1978), 483 - 497.

\bibitem{Oxt-Ul}
J.\,C.\,Oxtoby and S.\,M.\,Ulam, Measure preserving homeomorphisms and metrical transitivity, ``Ann. Math.'', \textbf{42} (1941), 874 - 920.

\bibitem{Yingst}
Andrew Q. Yingst, A characterization of homeomorphic Bernoulli trial measures, ``Trans. Amer. Math. Soc.'', \textbf{360} (2008), 1103 - 1131.

\bibitem{Yuasa}
Hisatoshi Yuasa, On the topological Orbit Equivalence in a Class of Substitution Minimal Systems, ``Tokyo J. Math.'', \textbf{25} (2002), 221 - 240.


\end{thebibliography}
\end{document}